\documentclass[final,onefignum,onetabnum]{siamart190516}

\usepackage[prependcaption,colorinlistoftodos]{todonotes}

\usepackage{cite}
\usepackage{amsmath,amsfonts,amssymb,amscd,latexsym,theorem,enumitem}
\usepackage{graphicx}
\usepackage{psfrag}
\usepackage{epsfig}
\usepackage{theorem}
\makeatletter
\xdef\@endgadget#1{{\unskip\nobreak\hfil\penalty50\hskip1em\hbox{}\nobreak
    \hfil#1\parfillskip=0pt\finalhyphendemerits=0\par}}
\def\@qedsymbol{${}_\blacksquare$}
\def\QED{\@endgadget{\@qedsymbol}}
\def\@Endofsymbol{$_\square$}
\def\Endoftheorem{\@endgadget{\@Endofsymbol}}

{\theorembodyfont{\itshape}}
{\theorembodyfont{\upshape}\newtheorem{remark}[theorem]{Remark}}
{\theorembodyfont{\upshape}\newtheorem{example}[theorem]{Example}}

\newcommand{\mR}{\mathbb{R}}

\newcommand{\pperp}{\perp \!\!\!\perp}
\newcommand{\X}{{\mathcal{X}}}
\newcommand{\dx}{\delta x}
\newcommand{\du}{\delta u}
\newcommand{\dy}{\delta y}
\newcommand{\pder}[2]{\frac{\partial #1}{\partial #2}}

\newcommand{\D}{{\mathcal{D}}}
\newcommand{\M}{{\mathcal{M}}}
\newcommand{\F}{{\mathcal{F}}}
\newcommand{\E}{{\mathcal{E}}}

\newcommand{\G}{{\mathcal{G}}}

\DeclareMathOperator{\im}{im}
\DeclareMathOperator{\col}{col}

\DeclareMathOperator{\rank}{rank}

\DeclareMathOperator{\clo}{cl}
\DeclareMathOperator{\rint}{rint}

\DeclareMathOperator{\dom}{dom}

\let\leq\leqslant
\let\geq\geqslant
\let\emptyset\varnothing

\newcommand{\calD}{\ensuremath{\mathcal{D}}}
\newcommand{\calE}{\ensuremath{\mathcal{E}}}
\newcommand{\calF}{\ensuremath{\mathcal{F}}}
\newcommand{\calG}{\ensuremath{\mathcal{G}}}
\newcommand{\calH}{\ensuremath{\mathcal{H}}}

\newcommand{\calM}{\ensuremath{\mathcal{M}}}
\newcommand{\calN}{\ensuremath{\mathcal{N}}}

\newcommand{\calR}{\ensuremath{\mathcal{R}}}
\newcommand{\calS}{\ensuremath{\mathcal{S}}}

\newcommand{\calV}{\ensuremath{\mathcal{V}}}

\newcommand{\bare}{\ensuremath{\bar{e}}}
\newcommand{\barf}{\ensuremath{\bar{f}}}

\newcommand{\bmat}{\begin{matrix}}
\newcommand{\emat}{\end{matrix}}
\newcommand{\bbm}{\begin{bmatrix}}
\newcommand{\ebm}{\end{bmatrix}}
\newcommand{\bpm}{\begin{pmatrix}}
\newcommand{\epm}{\end{pmatrix}}
\newcommand{\bse}{\begin{subequations}}
\newcommand{\ese}{\end{subequations}}
\newcommand{\beq}{\begin{equation}}
\newcommand{\eeq}{\end{equation}}
\newcommand{\ben}{\begin{enumerate}}
\newcommand{\een}{\end{enumerate}}
\newcommand{\beni}{\renewcommand{\labelenumi}{\roman{enumi}.}
\renewcommand{\theenumi}{\roman{enumi}}\begin{enumerate}}
\newcommand{\eeni}{\end{enumerate}\renewcommand{\labelenumi}{\arabic{enumi}.}
\renewcommand{\theenumi}{\arabic{enumi}}}
\newcommand{\bena}{\renewcommand{\labelenumi}{\alpha{enumi}.}
\renewcommand{\theenumi}{\alpha{enumi}}\begin{enumerate}}
\newcommand{\eena}{\end{enumerate}\renewcommand{\labelenumi}{\arabic{enumi}.}
\renewcommand{\theenumi}{\arabic{enumi}}}
\newcommand{\bit}{\begin{itemize}}
\newcommand{\eit}{\end{itemize}}
\newcommand{\bthe}{\begin{theorem}}
\newcommand{\ethe}{\end{theorem}}
\newcommand{\blem}{\begin{lemma}}
\newcommand{\elem}{\end{lemma}}
\newcommand{\bprop}{\begin{proposition}}
\newcommand{\eprop}{\end{proposition}}
\newcommand{\bex}{\begin{example}}
\newcommand{\eex}{\end{example}}
\newcommand{\bas}{\begin{assumption}}
\newcommand{\eas}{\end{assumption}}
\newcommand{\bre}{\begin{remark}}
\newcommand{\ere}{\end{remark}}
\newcommand{\bcor}{\begin{corollary}}
\newcommand{\ecor}{\end{corollary}}
\newcommand{\bdfn}{\begin{definition}}
\newcommand{\edfn}{\end{definition}}

\newcommand{\half}{\ensuremath{\frac{1}{2}}}

\newcommand{\pset}[1]{\ensuremath{\{#1\}}}

\newcommand{\set}[2]{\ensuremath{\{#1\mid #2\}}}

\newcommand{\abs}[1]{\ensuremath{| #1 |}}

\newcommand{\R}{\ensuremath{\mathbb R}}

\newcommand{\BP}{\noindent{\bf Proof. }}
\newcommand{\EP}{\hspace*{\fill} $\blacksquare$\bigskip\noindent}

\newcommand{\dualp}[2]{\ensuremath{\langle #1 \mid #2\rangle }}
\newcommand{\bilin}[2]{\ensuremath{\langle\langle #1,#2\rangle\rangle}}

\include{avds-newcommands}
\newcommand{\bq}{\begin{equation}}
\newcommand{\eq}{\end{equation}}
\newcommand{\bma}{\begin{bmatrix}}
\newcommand{\ema}{\end{bmatrix}}

\title{
Port-Hamiltonian systems and monotonicity
}

\author{M.K. Camlibel\thanks{Bernoulli Institute for Mathematics, Computer Science and Artificial Intelligence, University of Groningen, The Netherlands, (\email{m.k.camlibel@rug.nl})} \and A.J. van der Schaft\thanks{Bernoulli Institute for Mathematics, Computer Science and Artificial Intelligence, University of Groningen, The Netherlands}, \email{a.j.van.der.schaft@rug.nl}).}

\begin{document}
\maketitle
\thispagestyle{empty}
\pagestyle{empty}
\begin{abstract}
The relationships between port-Hamiltonian systems modeling and the notion of monotonicity are explored. The earlier introduced notion of incrementally port-Hamiltonian systems is extended to maximal cyclically monotone relations, together with their generating functions. This gives rise to new classes of incrementally port-Hamiltonian systems, with examples stemming from physical systems modeling as well as from convex optimization. An in-depth treatment is given of the composition of maximal monotone and maximal cyclically monotone relations, where in the latter case the resulting maximal cyclically monotone relation is shown to be computable through the use of generating functions.
Furthermore, connections are discussed with incremental versions of passivity, and it is shown how incrementally port-Hamiltonian systems with strictly convex Hamiltonians are (maximal) equilibrium independent passive. Finally, the results on compositionality of monotone relations are employed for a convex optimization approach to the computation of the equilibrium of interconnected incrementally port-Hamiltonian systems.   
\end{abstract}

\section{Introduction}
Port-based modeling of physical systems leads to their description as {\it port-Hamiltonian systems}. Such models have turned out to be powerful for purposes of analysis, simulation and control, see e.g. \cite{geoplexbook,passivitybook,vds-jeltsema}. On the other hand, during the last decades the concept of {\it monotonicity} has attracted much attention from multiple points of view. In relation with the current paper the following two aspects of monotonicity are most relevant. First, monotonicity has been a key notion in the study of nonlinear electrical circuits and general nonlinear network dynamics; see e.g. the recent paper \cite{chaffey} for a historical context and references. From a systems and control point of view this view on monotonicity is strongly related to notions of incremental passivity \cite{Desoer} and contraction \cite{lohmiller}. Second, monotonicity has evolved as a key concept in convex optimization (see e.g. \cite{parikh-boyd} and the references therein), as well as in nonlinear analysis (see e.g. \cite{brezis,bauschke}).

The present paper takes a closer look at the connections between port-Hamilto\-nian systems and monotonicity, and explores overarching notions. Already in our paper \cite{camvds}, inspired by \cite{j38}, we defined a new class of dynamical systems, coined as {\it incrementally port-Hamiltonian systems}. This was done by replacing the composition of the Dirac structure and the energy-dissipating relation in the standard definition of port-Hamiltonian systems by a general (maximal) monotone relation. Furthermore, in \cite{camvds} it was shown how monotone relations share the same compositionality property as Dirac structures, and sufficient conditions for the composition of two maximally monotone relations to be again maximally monotone were given. Moreover, the connections between incrementally port-Hamiltonian systems and the notions of incremental and differential passivity were briefly discussed. 
In the current paper this line of research is continued by developing a full-fledged theory of composition of maximally monotone relations. Furthermore, we take an in-depth look at (maximal) {\it cyclically} monotone relations in the context of port-Hamiltonian systems modeling. In particular, we show how, under mild technical conditions, the composition of two maximal cyclically monotone relations is again a maximal cyclically monotone relation. Maximal cyclically monotone relations are of special interest because they correspond to extended convex functions, in the sense that any maximal cyclically monotone relation is given as the subdifferential of a convex function, called the generating function of the relation. We show how the composition of two maximal cyclically monotone relations can be directly computed via their generating functions. Also, we define an appealing class of incrementally port-Hamiltonian systems which are defined by convex functions of the state and input. Obviously, this connects the theory of incrementally port-Hamiltonian systems to convex optimization. In fact, simple examples are continuous-time gradient algorithms for convex functions, and primal-dual gradient algorithms in case of minimization under affine constraints. Furthermore, we show how the equilibrium of interconnections of incrementally port-Hamiltonian systems defined by maximal cyclically monotone relations can be computed by convex optimization; thereby extending the innovative work of \cite{buergerzelazo}.
Finally, another connection with convex analysis appears if we assume the Hamiltonian function of the incrementally port-Hamiltonian system to be convex. This leads to {\it shifted passivity} \cite{passivitybook} of steady states, and in particular to (maximal) {\it equilibrium independent passivity} \cite{arcak, buergerzelazo}.

The organization of the paper is as follows. In Section~\ref{s:phs}, we quickly review the concepts of Dirac structures and standard port-Hamiltonian systems. This is followed by the definition of incrementally port-Hamiltonian systems in Section~\ref{s:inc-ph}, and the discussion of a number of examples. In Section~\ref{s:comp-mon}, we prove that under mild technical conditions the composition of two (maximal) (cyclically) monotone relations is (maximal) cyclically monotone, and thus the power-conserving interconnection of incrementally port-Hamiltonian systems is again an incrementally port-Hamiltonian system. In Section~\ref{s:ph-inc-ph} we investigate the various relationships and differences between incrementally port-Hamiltonian and standard port-Hamiltonian systems. As a direct application of the established compositionality theory we study in Section \ref{s:steady} the structure of the set of steady states of incrementally port-Hamiltonian systems, and the computation of the equilibrium of interconnected maximal cyclically monotone port-Hamiltonian systems via convex optimization. Furthermore, assuming convexity of the Hamiltonian, the relations with (maximal) equilibrium independent passivity are investigated. Finally in Section ~\ref{s:connections} the relations of incrementally port-Hamiltonian systems with incremental passivity and differential passivity are discussed, continuing upon the explorations in \cite{camvds}. 
The conclusions are in Section~\ref{s:conc}.

\section{Review of port-Hamiltonian systems on linear state spaces}\label{s:phs}
In order to motivate the definition of incrementally port-Hamiltonian systems we first review the definition of `ordinary' port-Hamiltonian systems; cf. \cite{vds-jeltsema, geoplexbook,passivitybook} for more details and ramifications.

Underlying the definition of a port-Hamiltonian system is the geometric notion of a {\it Dirac structure}, which relates the power variables of the constituting elements of the system in a power-conserving manner. Since incrementally port-Hamiltonian systems will be defined on {\it linear} state spaces we restrict as well attention to port-Hamiltonian systems on linear state spaces, and correspondingly to {\it constant} Dirac structures on linear spaces\footnote{For the extension to port-Hamiltonian systems on manifolds, and the corresponding notions of Dirac on manifolds we refer to e.g. \cite{Dalsmo, passivitybook, Courant}.}.
 
Power variables (such as voltages and currents, and forces and velocities), appear in conjugated pairs, whose products have dimension of power. They take values in dual linear spaces, with product meaning duality product. In particular, let $\F$ be a finite-dimensional linear space and $\mathcal{E} := \mathcal{F}^*$ be its dual space. We call $\mathcal{F}$ the space of {\it flow} variables, and $\mathcal{E}$ the space of {\it effort} variables. The duality product for the pair $(\E,\F)$, denoted by $\dualp{\cdot}{\cdot}$, is given as 
\[
\dualp{e}{f}=e^Tf\in\R
\]
for $e\in\E$ and $f\in\F$, and is the {\it power} associated to the pair $(f,e)$. Furthermore on $\mathcal{F} \times \mathcal{E}$ an {\it indefinite bilinear form} is defined as
$$
\bilin{(f_1,e_1)}{(f_2,e_2)} = \dualp{e_1}{f_2} + \dualp{e_2}{f_1},
$$
where $(f_i,e_i)\in\mathcal{F} \times \mathcal{E}$ with $i\in\pset{1,2}$. For any subspace $\calS \subset \mathcal{F} \times \mathcal{E}$, we denote its orthogonal companion with respect to this indefinite bilinear form by $\calS^{\pperp}$.

Throughout the paper, we will work with various spaces of flow/effort variables. By convention, if $\calF_\bullet$ denotes a certain space of flow variables then $\calE_\bullet:=\calF_\bullet^*$ will denote the corresponding space of effort variables.

\begin{definition}
Let $\mathcal{F}$ be a linear space. A subspace $\mathcal{D} \subset \mathcal{F} \times \mathcal{E}$ is a constant {\it Dirac structure} on $\mathcal{F}$ if $\mathcal{D} = \mathcal{D}^{\pperp}$.
\end{definition}
From now on in this paper a Dirac structure will simply refer to a {\it constant} Dirac structure on a linear space.
\begin{remark}\label{r:dirac}
An equivalent definition is the following \cite{Courant,geoplexbook,Dalsmo}. A Dirac structure is any subspace $\mathcal{D}$ with the property
\beq\label{e:dirac-prop}
\dualp{e}{f}=0 \mbox{ for all } (f,e) \in \cal D,
\eeq
which is {\it maximal} with respect to this property. (That is, there does not exist a subspace $\mathcal{D}'$ with $\mathcal{D} \subsetneq \mathcal{D}'$ such that $\dualp{e}{f}=0$ for all $(f,e) \in \mathcal{D}'$.)

In the finite-dimensional case (as will be the case throughout this paper) the maximal dimension of any subspace $\mathcal{D}$ satisfying \eqref{e:dirac-prop} equals $\dim \mathcal{F} = \dim \mathcal{E}$. Thus, equivalently, a Dirac structure is any subspace $\mathcal{D}$ satisfying \eqref{e:dirac-prop} together with 
\[
\dim \mathcal{D} = \dim \mathcal{F}.
\] 
\end{remark}

The definition of a {\it port-Hamiltonian system} on a linear space contains the following ingredients (see e.g. \cite{vds-maschke,passivitybook,Dalsmo,vds-jeltsema}). First a Dirac structure $\mathcal{D}$ defined on the space of all flow variables, that is,
\begin{equation}\label{diracres}
\mathcal{D} \subset \mathcal{F}_x \times  \mathcal{F}_P \times \mathcal{F}_R \times \mathcal{E}_x \times \mathcal{E}_P \times \mathcal{E}_R
\end{equation}
Here $(f_x,e_x) \in \mathcal{F}_x \times \mathcal{E}_x$ are the flow and effort variables linking to the {\it energy-storing} elements, $(f_R,e_R) \in \mathcal{F}_R \times \mathcal{E}_R$ are the flow and effort variables linking to {\it energy-dissipating} elements, and finally $(f_P,e_P) \in \mathcal{F}_P \times \mathcal{E}_P$ are the flow and effort {\it port} variables.
The port-Hamiltonian system is defined by specifying, next to its Dirac structure $\mathcal{D}$, the constitutive relations of the energy-dissipating elements, and of the energy-storing elements. An {\it energy-dissipating relation} is any subset $\mathcal{R} \subset \mathcal{F}_R \times \mathcal{E}_R$ with the property
\begin{equation}\label{res}
\dualp{e_R}{f_R} \geq 0 \mbox{ for all } (f_R,e_R) \in \mathcal{R}.
\end{equation}
The constitutive relations of the {\it energy-storing} elements are specified by a {\it Hamiltonian} $H: \mathcal{X} \to \mathbb{R}$, where $\mathcal{X}=\mathcal{F}_x$.
Thus the total energy while at state $x$ is given as $H(x)$. This defines the following constitutive relations between the state variables $x$ and the flow and effort vectors $(f_x,e_x)$ of the energy-storing elements\footnote{Throughout this paper the vector $\frac{\partial H}{\partial x}(x)$ denotes the {\it column} vector of partial derivatives; the corresponding row vector denoted as $\frac{\partial H}{\partial x^T}(x)$.}
\begin{equation}
\label{storage}
\dot{x} = -f_x \quad\text{and}\quad e_x = \frac{\partial H}{\partial x}(x).
\end{equation}

\begin{definition}
Consider a Dirac structure (\ref{diracres}), a Hamiltonian $H:\mathcal{X} \to \mathbb{R}$, and an energy-dissipating relation $\mathcal{R} \subset \mathcal{F}_R \times \mathcal{E}_R$ as above. Then the dynamics of the corresponding {\it port-Hamiltonian system} on $\mathcal{X}$ is given as
\bse
\label{switchingPH}
\begin{gather}
\left(-\dot{x}(t),f_P(t),-f_R(t),\frac{\partial H}{\partial x}\big(x(t)\big),e_P(t),e_R(t)\right) \in \D, \\ \big(f_R(t),e_R(t) \big) \in \mathcal{R}
\end{gather}
\ese
at (almost) all time instants $t$. 
\end{definition}
Equation \eqref{storage} immediately implies the energy balance
$\frac{d}{dt}H = \frac{\partial H}{\partial x^T}(x) \dot{x} = -\dualp{e_x}{f_x}$. 
Furthermore, the {\it composition} 
\begin{gather}
\mathcal{D} \rightleftarrows \mathcal{R} := \{ (f_x, f_P, e_x, e_P) \in  \mathcal{F}_x\times \mathcal{F}_P \times \mathcal{E}_x \times\mathcal{E}_P\mid \notag\\[2mm]
\!\!\!\!\!\!\!\exists (f_R, e_R) \in \mathcal{R}  
\mbox{ s.t. }  (f_x, f_P, -f_R, e_x, e_P,e_R) \in \mathcal{D} \}
 \end{gather}
 satisfies by the power-conserving property of the Dirac structure and \eqref{res}
 \begin{equation}
 e_x^Tf_x + e_P^Tf_P = e_R^Tf_R \geq 0 
  \end{equation}
for all $(f_x, f_P, e_x, e_P) \in \mathcal{D} \rightleftarrows \mathcal{R}$. Taken together this implies that 
\begin{equation}\label{energybalance2}
\frac{d}{dt} H(x(t)) \leq e_P^T(t)f_P(t),
\end{equation}
showing {\it cyclo-passivity} of any port-Hamiltonian system, and {\it passivity}  if $H:\mathcal{X} \to \mathbb{R}_+$ \cite{passivitybook} .

\section{Incrementally port-Hamiltonian systems}\label{s:inc-ph}
The basic idea in the definition of an {\it incrementally port-Hamiltonian system}, as first introduced in \cite{camvds}, is to replace the composition $\mathcal{D} \rightleftarrows \mathcal{R}$ of a Dirac structure $\mathcal{D}$ and an energy-dissipating relation $\mathcal{R}$ by a {\it monotone relation} $\M$. To do so, we begin with a quick review of monotone relations.

\begin{definition}\label{d:monotone}
A relation $\M\subset \F\times \E$ is said to be 
\begin{itemize}
\item {\em monotone\/} if
\[
\dualp{e_1-e_2}{f_1-f_2}\geq 0
\]
for all $(f_i,e_i)\in \M$ with $i\in\pset{1,2}$.
\item {\em cyclically monotone\/} if 
$$
\dualp{e_0}{f_0-f_1}+\dualp{e_1}{f_1-f_2}+\cdots+\dualp{e_{m-1}}{f_{m-1}-f_m}+\dualp{e_m}{f_m-f_0}\geq 0.
$$
for all $m\geq 1$ and $(f_i,e_i)\in \M$ with $i\in\pset{0,1,\ldots,m}$.
\end{itemize}
\end{definition}
Since $
\dualp{e_0}{f_0-f_1}+\dualp{e_1}{f_1-f_0}=\dualp{e_0-e_1}{f_0-f_1}$
for all $e_0,f_0,e_1,f_1$, every cyclically monotone relations is automatically monotone. 

A simple example of a monotone relation $\mathcal{M} \subset \mathbb{R}\times\R$ is the graph of a monotone (i.e., non-decreasing), possibly discontinuous, function. For example, the graph of the discontinuous function $\theta: \mathbb{R} \to \mathbb{R}$ given by
\begin{equation}\label{eq:disc}
\theta(x) = \begin{cases}
-1&\text{ if }x<0\\
1 & \text{ if }x\geq 0
\end{cases}
\end{equation}
is a monotone relation. This example already motivates the strengthened definition of a {\it maximal} monotone relation.
\begin{definition}
A relation $\M\subset \F\times \E$ is called {\em maximal (cyclically) monotone\/} if it is (cyclically) monotone and the implication
\[
\M'\text{ is (cyclically) monotone and } \M \subset \M'\quad \implies\quad \M=\M'
\]
holds. 
\end{definition}
The graph of the discontinuous function $\theta$ in (\ref{eq:disc}) is monotone, but not maximal monotone. In fact, its graph can be enlarged so as to obtain the following maximal monotone relation
\begin{equation}\label{eq:maxmon}
\mathcal{M} = \left\{ (x,y)\mid
y\in\begin{cases}
\pset{-1}&\text{ if }x<0\\
[-1,1] & \text{ if }x= 0\\
\pset{1}&\text{ if }x>0
\end{cases}\right\}.
\end{equation}
Note that the function $\theta$ in (\ref{eq:disc}) can be regarded as the description of a relay, while its closure given by the maximal monotone relation $\mathcal{M}$ defined in (\ref{eq:maxmon}) defines for example an ideal Coulomb friction characteristic. 

A few well-known facts are noteworthy. For {\it continuous} functions, monotonicity of the graph implies maximal monotonicity (see e.g. \cite{rockafellar}). Also, every maximal monotone relation on $\R\times\R$ is maximal cyclically monotone. (Hence the above Coulomb friction characteristic in \eqref{eq:maxmon} is maximal cyclically monotone.)
In higher dimensions, however, not every maximal monotone relation enjoys the cyclical monotonicity property. Indeed, for example the relation given by
\[
\left\{\left(\bbm x\\y\ebm,\bbm-y\\x\ebm\right)\mid
x,y\in\R
\right\}\subset\R^2\times\R^2
\]
is maximal monotone but not cyclically monotone. More generally, later on (Proposition~\ref{l:every dirac}) we will see that Dirac structures are maximal monotone, but not cyclically monotone if they are the graph of a non-zero map.

The importance of maximal cyclically monotone relations $\M$ lies in the fact that they correspond to extended real-valued {\it convex functions}. This will be briefly review next, for more details we refer to \cite{rockafellar}. Let $\phi:\F\rightarrow(-\infty,+\infty]$ be a proper convex function. Its {\em effective domain\/} is defined by
\[
\dom \phi:=\set{f\in\F}{\phi(f)<+\infty},
\]
its {\em subdifferential of $\phi$ at $f$\/} by
\[
\partial \phi(f):=\begin{cases}
\set{e\in\E}{\phi(\barf)\geq \phi(f)+\dualp{e}{\barf-f}\,\,\,\forall\,\barf\in\F}&\text{ if }f\in\dom(\phi)\\
\emptyset&\text{ otherwise},
\end{cases}
\]
and its {\em conjugate\/} $\phi^\star:\E\rightarrow(-\infty,+\infty]$ by
\[
\phi^\star(e):=\sup\set{\dualp{e}{f}-\phi(f)}{f\in\F}.
\]
If, in addition, $\phi$ is lower semicontinuous, then $\phi=(\phi^{\star})^{\star}$ and 
\beq\label{e:inverse partial duality}
e\in\partial \phi(f) \iff f\in\partial \phi^\star(e).
\eeq
It turns out (see \cite[Thm. 12.25]{rockafellar}) that a relation $\calM\subset\F\times\E$ is maximal cyclically monotone if and only if there exists a proper lower semicontinuous convex function $\phi$ such that 
\[
\calM=\set{(f,e)}{e\in\partial \phi(f)}=\set{(f,e)}{f\in\partial \phi^\star(e)}.
\]
In this case, we say that $\calM$ is {\em generated by\/} $\phi$, or that $\phi$ is a {\em generating function} of $\calM$. Note that $\phi$ is determined by $\calM$ uniquely up to an additive constant. 

As an example, consider the relation given by \eqref{eq:maxmon}. One easily verifies that $\calM$ is generated by the convex function $\phi(x)$ given by $x \mapsto\abs{x}$. Furthermore, 
\[\phi^\star(y)=\begin{cases}0&\text{if }y\in[-1,1]\\
+\infty&\text{if }y\not\in[-1,1].\end{cases}\]
The definition of an incrementally port-Hamiltonian system as introduced in \cite{camvds} is now extended as follows.
\begin{definition}
Consider a (maximal) (cyclically) monotone relation
\[
\mathcal{M} \subset \mathcal{F}_x \times \mathcal{E}_x \times \mathcal{F}_P \times \mathcal{E}_P
\]
and a Hamiltonian $H:\mathcal{F}_x \to \mathbb{R}$. Then the dynamics of the corresponding {\it (maximal) (cyclically) monotone port-Hamiltonian system}, is given by the requirement
\begin{equation}\label{inc-PH}
\big(-\dot{x}(t),f_P(t),\frac{\partial H}{\partial x}\big(x(t)\big),e_P(t)\big) \in \M \\
\end{equation}
for all time instants $t$.
\end{definition}
\begin{remark}
Throughout the terminology {\it incrementally port-Hamiltonian system} will be used as shorthand for all systems defined with respect to monotone relations $\M$. Whenever we need to be more precise about the properties of the monotone relation $\M$ we will refer to the system as a (maximal) (cyclically) monotone port-Hamiltonian system.
\end{remark}
It follows that the dynamics of any incrementally port-Hamiltonian system satisfies the {\it incremental dissipation inequality}
\begin{equation}\label{e:disp-like}
\dualp{\frac{\partial H}{\partial x}\big(x_1\big)-\frac{\partial H}{\partial x}\big(x_2\big)} {\dot{x}_1-\dot{x}_2} \leq \dualp{e_P^1-e_P^2}{f_P^1-f_P^2}
\end{equation}
for all quadruples $(x_i, \dot{x}_i, f_P^i,e_P^i), i=1,2,$ satisfying
$
\big(-\dot{x}_i,f_{P}^i,,\frac{\partial H}{\partial x}\big(x_i\big),e_{P}^i\big) \in \calM
$
with $i\in\pset{1,2}$. The consequences of this dynamical inequality, and especially the relation with incremental and differential passivity, will be discussed in Section \ref{s:connections}.

Incrementally port-Hamiltonian systems are ubiquitous in physical systems modeling, as already illustrated by the following examples.
\begin{example}[Mechanical systems with friction]
Consider a mechanical system subject to friction. The friction characteristic is given by a relation between $f_R$, $e_R$. In the case of a scalar friction characteristic of the form $e_R = R(f_R)$ the system is port-Hamiltonian if the graph of the function $R: \mR \to \mR$ is in the first and third quadrant. On the other hand, it is maximal cyclically monotone port-Hamiltonian if the function $R$ is a monotonically non-decreasing and moreover continuous, or otherwise the graph of $R$ is extended by the interval between the left- and right limit values at its discontinuities. (A typical example of the latter is {\it Coulomb friction} as mentioned above.)
\end{example}

\begin{example}[Systems with constant sources]
Physical systems containing nonzero internal {\it constant sources} are {\it not} port-Hamiltonian but can be incrementally port-Hamiltonian. Consider for example any LC-circuit with passive resistors/conductors and constant voltage and/or current sources. The same holds for an arbitrary mechanical system with constant actuation: incrementally port-Hamiltonian but {\it not} port-Hamiltonian for nonzero constant actuation.
\end{example}

\begin{example}[Van der Pol oscillator]
Consider an electrical LC-circuit (with possibly nonlinear capacitors and inductors), together with a single conductor with current $f_R=I$ and voltage $e_R=V$. In case of a linear conductor $I=GV, G>0$, the system is both port-Hamiltonian and maximal monotone port-Hamiltonian. For a nonlinear conductor $I = G(V)$ the system is port-Hamiltonian if and only the graph of the function $G$ is in the first and third quadrant and maximal monotone port-Hamiltonian if $G$ is monotonically non-decreasing and continuous, or otherwise the graph of $G$ is extended by the interval between the left- and right limit values at its discontinuities.
For example, the conductor characteristic $I=\Phi(V)$ where $\Phi(z)= \gamma z^3 - \alpha z, \alpha, \gamma >0,$ defines a system which is port-Hamiltonian but {\it not} monotone port-Hamiltonian, since the function $\Phi$ is not monotone. On the other hand, by adding a constant source voltage $V_0$ and constant source current $I_0$ in such a way that the {\it tunnel diode} characteristic
\[
I = \Phi (V - V_0) + I_0,
\]
passes through the origin the resulting system (the Van der Pol oscillator) is {\it not} port-Hamiltonian, since close to the origin the characteristic is in the second and fourth quadrant, while neither is it incrementally port-Hamiltonian. 
\end{example}

An appealing class of maximal cyclically monotone port-Hamiltonian systems is defined as follows. Consider any Hamiltonian $H: \mathcal{X} \to \mathbb{R}$, and any {\it convex} function $K: \mathcal{X}^* \times \mathcal{U} \to \mathbb{R}$. Then the system
\begin{equation} 
\label{systemconvex}
\begin{array}{rcl}
\dot{x} & = & -\frac{\partial K}{\partial e} (\nabla H(x),u), \quad e = \nabla H(x):=\frac{\partial H}{\partial x}(x) \\[2mm]
y & = & \frac{\partial K}{\partial u} (\nabla H(x),u)
\end{array}
\end{equation}
is a maximal cyclically monotone port-Hamiltonian system with maximal cyclically monotone relation
$\mathcal{M} = \mbox{graph } (\partial K)$. Special case occurs if the convex function $K(e,u)$ is of the form
\begin{equation}
K(e,u) = P(e) + e^TBu,
\end{equation}
with $P$ a convex function of $e$, and $B$ an $n \times m$ matrix. This yields the restricted system class
\begin{equation} 
\begin{array}{rcl}
\dot{x} & = & -\frac{\partial P}{\partial e} (\nabla H(x)) - Bu\\[2mm]
y & = & B^T \nabla H(x)
\end{array}
\end{equation}

A physical example of the form \eqref{systemconvex} is the following.
\begin{example}[Nonlinear RC electrical circuit]
Consider an RC electrical circuit, with nonlinear conductors at the edges and grounded nonlinear capacitors at part of the nodes, while the remaining nodes are the boundary nodes (terminals). Let the circuit graph be defined by an incidence matrix $D$, split according to the splitting of the capacitor and boundary nodes as
\bq
D = \bma D_c \\ D_b \ema
\eq
Furthermore, let the conductors at the edges be given as $I_j = G_j(V_j)$, where $I_j,V_j$ are the current through, respectively, voltage, across the $j$-th edge, $j=1, \cdots,m$. Assume that the conductance functions $G_j$ are all monotone (however not necessarily in the first and third quadrant). This means that there exist convex functions $\hat{K}_j$ such that $G_j(V_j)= \frac{d\hat{K}_j}{dV_j}(V_j)$ (if for simplicity we assume that the functions $G_j$ are continuous, and $\hat{K}_j$ are differentiable). Define the convex functions
\bq
\hat{K}(V_1, \ldots, V_m) := \sum_{j=1}^m \hat{K}_j(V_j), \quad K(\psi) := \hat{K}(D^T \psi),
\eq
where $\psi$ is the vector of node voltage potentials. (Recall that by Kirchhoff's voltage law $V=D^T \psi$.) It is immediately checked that $\frac{\partial K}{\partial \psi} = D \frac{\partial \hat{K}}{\partial V}(D^T \psi)$.
Denote the vector of charges of the grounded capacitors by $Q$. It follows by Kirchhoff's current laws that the dynamics of the nonlinear RC circuit is given by
\bq
\label{eq:RC}
\begin{array}{rcl}
\dot{Q} & = & - D_c \frac{\partial \hat{K}}{\partial V}(D^T \psi) \\[2mm]
I_e & = & D_e \frac{\partial \hat{K}}{\partial V}(D^T \psi) ,
\end{array}
\eq
where $I_e$ is the vector of injected currents at the boundary nodes of the electrical circuit. According to the splitting of the nodes in capacitor and boundary nodes write $\psi = \bma \psi_c \\ \psi_e \ema$. Then by specifying the nonlinear grounded capacitors by a Hamiltonian function $H(Q)$ it follows that $\psi_c =\frac{\partial H}{\partial Q}(Q)$.

The system \eqref{eq:RC} is a maximal cyclically monotone port-Hamiltonian system of the form \eqref{systemconvex}, with inputs $\psi_e$, state $Q$ (dimension equal to the number of capacitor nodes), and outputs $I_e$. The generating function of its maximal cyclically monotone relation is given by the convex function $K(\psi)$. Finally note that the system equations can be also written in terms of the alternative state vector $\psi_c$ (under the assumption that the map $Q \mapsto \frac{\partial H}{\partial Q}(Q)$ is invertible), by substituting $\dot{\psi}_c = \frac{\partial^2 H}{\partial Q^2}(Q) \dot{Q}$.
\end{example}

An example of a maximal cyclically monotone port-Hamiltonian system of the form \eqref{systemconvex} that is {\it not} stemming from physical systems modeling, but instead from {\it optimization}, is the following.

\begin{example}[Gradient algorithm in continuous time]
Consider the problem of minimizing a convex function $P: \mR^n \to \mR$. The gradient algorithm in continuous time is given as
\beq
\begin{array}{rclrl}
\tau \dot{q}  & = & - \frac{\partial P}{\partial q}(q)  - Bu \\[2mm]
y & = & B^T q \, , &&
\end{array}
\eeq
where $\tau$ is a positive definite matrix determining the time-scales of the algorithm.
Here, an input vector $u \in \mR^n$ is added in order to represent possible interaction with other algorithms or dynamics (e.g., if the gradient algorithm is carried out in a distributed fashion), defining a conjugate output vector as $y = B^T q \in \mR^n$.

This defines a maximal cyclically monotone port-Hamiltonian system with state vector $x:= \tau q$, quadratic Hamiltonian $H(x)= \frac{1}{2}x^T  \tau^{-1} x$, and  maximal cyclically monotone relation
\beq
\M = \{ (f_S, e_S, y,u) \mid -f_S= \frac{\partial P}{\partial q}(q) + Bu, y= B^T e_S\},
\eeq
where $e_S = \frac{\partial H}{\partial x}(x)= \tau^{-1} x$.
\end{example}

An {\it extended} class of maximal monotone port-Hamiltonian systems is defined as
\begin{equation} 
\begin{array}{rcl}
\dot{x} & = & J\nabla H(x) - \frac{\partial P}{\partial e} (\nabla H(x)) - Bu\\[2mm]
y & = & B^T \nabla H(x)
\end{array}
\end{equation}
where $J$ is a skew-symmetric matrix, and $P$ a convex function as above. It will follow from Proposition \ref{l:every dirac} that if $J \neq 0$ then the underlying maximal monotone relation is {\it} not derivable from a convex function, and the system is {\it not} cyclically monotone port-Hamiltonian anymore. An example within this class is the following.

\begin{example}[Primal-dual gradient algorithm \cite{arrow}]
Consider the {\it constrained} optimization problem 
\beq
\min_{q; \, Aq=b} P(q) ,
\eeq
where $P: \mR^n \to \mR$ is a convex function, and $Aq=b$ are affine constraints for some $k \times n$ matrix $A$ and vector $b \in \mR^k$. The resulting Lagrangian function is defined as
\beq
L(q,\lambda):= P(q) + \lambda^T(Aq - b), \quad \lambda \in \mR^k,
\eeq
which is convex in $q$ and concave in $\lambda$. The {\it primal-dual gradient} algorithm for solving the optimization problem in continuous time is given as
\beq
\begin{array}{rclrl}
\tau_q \dot{q} & = & - \frac{\partial L}{\partial q}(q,\lambda) & = & - \frac{\partial P}{\partial q}(q) - A^T \lambda + u \\[2mm]
\tau_{\lambda} \dot{\lambda} & = & \frac{\partial L}{\partial \lambda}(q,\lambda) &= & Aq - b \\[2mm]
y & = & q \, , &&
\end{array}
\eeq
where $\tau_q, \tau_{\lambda}$ are positive-definite matrices determining the time-scales of the algorithm.
Again, an input vector $Bu \in \mR^n$ is added in order to represent possible interaction with other algorithms or dynamics, defining a conjugated output vector $y = G^T q \in \mR^n$.

This defines a maximal monotone port-Hamiltonian system with state vector $x = (x_q,x_{\lambda}):=(\tau_q q, \tau_{\lambda} \lambda)$, quadratic Hamiltonian 
\beq
H(x)= \frac{1}{2}x_q^T \tau_q^{-1} x_q +  \frac{1}{2} x_{\lambda} \tau^{-1}_{\lambda} x_{\lambda} ,
\eeq
and maximal monotone relation
\beq
\M = \{ (f_S, e_S, y,u) \mid -f_S= \begin{bmatrix} 0 & A^T \\[2mm] -A & 0 \end{bmatrix} 
e_S - \begin{bmatrix} \frac{\partial P}{\partial q}(q) \\[2mm] b \end{bmatrix} + \begin{bmatrix}0 \\[2mm] B \end{bmatrix}u, y=  \begin{bmatrix}0 & B^T \end{bmatrix}e_S\},
\eeq
where 
\bq
e_S = \nabla H(x)= \bma \tau^{-1}_q x_q \\ \tau^{-1}_{\lambda} x_{\lambda} \ema = \bma q \\ \lambda \ema
\eq
See for an application the optimization of social welfare in a dynamic pricing algorithm for power networks \cite{stegink}.
\end{example}

\section{Composition of monotone relations}\label{s:comp-mon}
A cornerstone of port-Hamiltonian systems theory is the fact that the power-conserving interconnection of port-Hamilto\-nian systems defines again a port-Hamiltonian system. This in turn is based on the fact that the {\it composition} of Dirac structures is again a Dirac structure. In this section we will show that the same property holds for incrementally port-Hamiltonian systems. This follows from the corresponding compositionality property of (maximal) (cyclically) monotone relations.

Let us start by considering two monotone relations $\calM_a \subset \F_a \times \F \times \E_a \times \E$ and $\calM_b \subset \F_b \times \F\times \E_b \times \E$. Define the {\it composition} of $\calM_a$ and $\calM_b$, denoted as $\calM_a \overset{\scriptscriptstyle\calF\times\calE}{\leftrightarrows}\calM_b$, as before, by
\begin{align}
&\calM_a \overset{\scriptscriptstyle\calF\times\calE}{\rightleftarrows} \calM_b  :=  \{(f_a,f_b,e_a,e_b) \in \F_a \times \F_b \times \E_a \times \E_b \mid \notag\\
&\exists (f,e)\in \F\times\E \text{ s.t. } 
(f_a,f,e_a,e) \in \calM_a ,
(f_b,-f,e_b,e) \in \calM_b \}.
\end{align}
Thus the composition of $\calM_a$ and $\calM_b$ is obtained by imposing the interconnection constraints
\begin{equation}\label{int}
f_1 = - f_2, \quad e_1=e_2,
\end{equation}
on the vectors $(f_a,f_1,e_a,e_1) \in \calM_a$ and $(f_b,f_2,e_b,e_2)  \in \calM_b$  
and looking at the resulting vectors $(f_a,f_b,e_a,e_b) \in \F_a \times \F_b \times \E_a \times \E_b$. 

Whenever interconnection flow and effort spaces $\calF$ and $\calE$ are clear from the context, we will simply write $\calM_a \leftrightarrows\calM_b$. The following result is straightforward.

\begin{proposition}
Let $\calM_a \subset \F_a \times \F \times \E_a \times \E$ and $\calM_b \subset \F_b \times \F\times \E_b \times \E$ be (cyclically) monotone relations. Then, $\calM_a \rightleftarrows \calM_b\subset\F_a \times \F_b \times \E_a \times \E_b$ is (cyclically) monotone.
\end{proposition}
\BP Suppose that both $\calM_a$ and $\calM_b$ are monotone relations. Let 
\[(f_a,f_b,e_a,e_b),(\barf_a,\barf_b,\bare_a,\bare_b)\in\calM_a \rightleftarrows \calM_b.
\]
Then, there exist $(f,e),(\barf,\bare)\in\F\times\E$ such that $(f_a,f,e_a,e),(\barf_a,\barf,\bare_a,\bare) \in \calM_a$ and $(f_b,-f,e_b,e),(\barf_b,-\barf,\bare_b,\bare) \in \calM_b$. From monotonicity of $\calM_a$ and $\calM_b$, we have
\[
\dualp{\bbm e_a-\bare_a\\e-\bare\ebm}{\bbm f_a-\barf_a\\f-\barf\ebm}\geq 0\quad\text{and}\quad
\dualp{\bbm e_b-\bare_b\\e-\bare\ebm}{\bbm f_b-\barf_b\\-f+\barf\ebm}\geq 0.
\]
By adding these left hand sides of these inequalities, we obtain
\[
\dualp{\bbm e_a-\bare_a\\e_b-\bare_b\ebm}{\bbm f_a-\barf_a\\f_b-\barf_b\ebm}\geq 0.
\]
This means that $\calM_a \rightleftarrows \calM_b$ is monotone. The cyclical monotone case follows in a similar fashion.\EP

Also the composition of two {\it maximal} monotone relations turns out to be maximal monotone; provided certain (mild) regularity conditions are met. To elaborate on this, we first introduce some nomenclature and review some known facts about maximal monotone relations.

For a set $S\in\F$, $\clo S$ denotes its {\em closure}. The {\em relative interior\/} of a convex set $C\subseteq\F$ is denoted by $\rint C$. A set $S\subseteq\F$ is said to be {\em nearly convex\/} if there exists a convex set $C\subseteq\F$ such that $C\subseteq S\subseteq \clo C$. For a nearly convex set $S$, in general, there can be multiple convex sets $C$ satisfying $C\subseteq S\subseteq \clo C$. For any such set $C$, however, we have that $\clo C=\clo S$. As such, $\clo S$ is convex if $S$ is nearly convex. Based on this observation, 
one can extend the notion of relative interior to nearly convex sets by defining $\rint S=\rint(\clo S)$.

Let $S\subseteq\F_1\times\F_2\times\E_1\times\E_2$. The projection of $S$ on $\F_1\times\F_2$, denoted by $\Pi(S,\F_1\times\F_2)$, is defined as
\[
\Pi(S,\F_1\times\F_2):=\set{(f_1,f_2)}{\exists\,(e_1,e_2)\in\E_1\times\E_2\text{ s.t. }(f_1,f_2,e_1,e_2)\in S}.
\]
We define projections of $S$ on $\E_1\times\E_2$, $\F_i\times\E_j$, $\F_i$, and $\E_j$ in a similar fashion.

Let $S\subseteq\calF\times\calG$ be a nearly convex set. Then, both $\Pi(S,\calF)$ and $\Pi(S,\calG)$ are nearly convex sets. Furthermore, one can show that 
\beq\label{e:nconvex relint}
\rint S=\set{(f,g)}{f\in\rint\Pi(S,\calF)\text{ and }g\in\rint\Pi(S\cap(\pset{f}\times\calG),\calG)}.
\eeq

Let $\calM\subset\F\times\E$ be a maximal monotone relation. Then, the projections $\Pi(\calM,\F)$ and $\Pi(\calM,\E)$ are nearly convex sets \cite[Thm. 12.41]{rockafellar}.

Let $L:\calG\rightarrow\calH$ be a linear map and $L^*:\calH^*\rightarrow\calG^*$ denote its adjoint. For maximal monotone relations $\calM\subseteq\calH\times\calH^*$ and $\calN\subseteq\calG\times\calG^*$,
define $\calM_L\subseteq\calG\times\calG^*$ and $_{L}\calN\subseteq\calH\times\calH^*$ by
\begin{gather*}
\calM_L=\set{(g,L^*h^*)}{(Lg,h^*)\in\calM}\\
_{L}\calN=\set{(Lg,h^*)}{(g,L^*h^*)\in\calN}
\end{gather*}
From \cite[Thm. 12.43]{rockafellar}, we know that $\calM_L$ is maximal monotone if
\beq\label{e:im L domain}
\im L\cap\rint\Pi(\calM,\calH)\neq\emptyset
\eeq
and $_{L}\calN$ is maximal monotone if
\beq\label{e:im L* domain}
\im L^*\cap\rint\Pi(\calN,\calG^*)\neq\emptyset.
\eeq

Furthermore, if $\calM$ is generated by $\phi:\calH\rightarrow(-\infty,+\infty]$ and
\beq\label{e:im L domain-subdiff}
\im L\cap\rint\dom(\phi)\neq\emptyset,
\eeq
then $\calM_L$ is generated by $\phi\circ L$ given by $h\mapsto \phi(Lh)$. Dually, if $\calN$ is generated by $\psi:\calG\rightarrow(-\infty,+\infty]$ and
\beq\label{e:im L* domain-subdiff}
\im L^*\cap\rint\dom(\psi^\star)\neq\emptyset,
\eeq
then $_L\calN$ is generated by the function $(\psi^\star\circ L^*)^\star$.

\bthe\label{t:mm preserve} 
Let $\calM_a \subset \F_a \times \F \times \E_a \times \E$ and $\calM_b \subset \F_b \times \F\times \E_b \times \E$ be maximal monotone relations. Let
\[
C_f=\set{(f_1,f_2)}{f_1\in\Pi(\calM_a,\F)\text{ and }f_2\in\Pi(\calM_b,\F)}.
\]
and
\[
C_e=\set{(e_1,e_2)}{\exists\,f\text{ s.t. } (f,e_1)\in\Pi(\calM_a,\F\times\E)\text{ and }(-f,e_2)\in\Pi(\calM_b,\F\times\E)}.
\]
Suppose that there exists $(\barf,\bare)\in\F\times\E$ such that
\ben[label=(\roman*), ref=(\roman*)]
\item\label{i:mm cond 1} $(\barf,-\barf)\in\rint C_f$ and
\item\label{i:mm cond 2} $(\bare,\bare)\in\rint C_e$.
\een
Then, $\calM_a \rightleftarrows \calM_b\subset\F_a \times \F_b \times \E_a \times \E_b$ is a maximal monotone relation. 
\ethe
\BP
First, we give an alternative characterization of $\calM_a\rightleftarrows\calM_b$. Let $\calM \subset \F_a \times \F \times \F_b \times \F\times \E_a \times \E\times \E_b \times \E$ be defined by
\[
\calM:=\set{(f_a,f_1,f_b,f_2,e_a,e_1,e_b,e_2)}{(f_a,f_1,e_a,e_1)\in\calM_a\text{ and }(f_b,f_2,e_b,e_2)\in\calM_b}.
\]
Since $\calM_a$ and $\calM_b$ are both maximal monotone, so is $\calM$. Let $A:\F_a\times\F\times\F_b\rightarrow\F_a\times\F\times\F_b\times\F$ be the linear map given by
\[
(f_a,f,f_b)\mapsto(f_a,f,f_b,-f)
\]
and $B:\F_a\times\F\times\F_b\rightarrow\F_a\times\F_b$ be the linear map given by
\[
(f_a,f,f_b)\mapsto (f_a,f_b).
\]
Note that $A^*:\E_a\times\E\times\E_b\times\E\rightarrow \E_a\times\E\times\E_b$ is given by
\[
(e_a,e_1,e_b,e_2)\mapsto(e_a,e_1-e_2,e_b) 
\]
and $B^*:\E_a\times\E_b\rightarrow\E_a\times\E\times\E_b$ given by
\[
(e_a,e_b)\mapsto(e_a,0,e_b).
\]
Now, we claim that
\[
\calM_a\rightleftarrows\calM_b=_{B}(\calM_A)
\]
To see this, note that
\begin{align*}
\calM_A&=\set{(f_a,f,f_b,e_a,e_1-e_2,e_b)}{(f_a,f,f_b,-f,e_a,e_1,e_b,e_2)\in\calM}\\
&=\{(f_a,f,f_b,e_a,e_1-e_2,e_b)\mid(f_a,f,e_a,e_1)\in\calM_a\\
&\qquad\text{ and }(f_b,-f,e_b,e_2)\in\calM_b\}
\end{align*}
and
\begin{align}
_{B}(\calM_A)&=\set{(f_a,f_b,e_a,e_b)}{\exists\,f\in\F\text{ s.t. }(f_a,f,f_b,e_a,0,e_b)\in\calM_a}\notag\\
&=\{(f_a,f_b,e_a,e_b)\mid\exists\,(f,e)\in\F\times\E\text{ s.t. }(f_a,f,e_a,e)\in\calM_a\notag\\
&\qquad\text{ and }(f_b,-f,e_b,e)\in\calM_b\}\notag\\
&=\calM_a\rightleftarrows\calM_b\label{e:BMA is Ma comp Mb}.
\end{align}
where $A:\calF_a\times\calV_a\times\calF_b\times\calV_b\rightarrow\calF_a\times\calV_a\times\calV_c\times\calV_d\times\calF_b\times\calV_b$ is the linear map given by
\[
(f_a,v_a,f_b,v_b)\mapsto(f_a,v_a,v_a,v_b,f_b,v_b)
\]
Since $\calM$ is maximal monotone, we see from \eqref{e:im L domain} that $\calM_A$ is maximal monotone if
\beq\label{e:maxmon cond for MA} 
\im A\cap\rint \Pi(\calM,\F_a\times\F\times\F_b\times\F)\neq\emptyset.
\eeq
From \eqref{e:nconvex relint}, it follows that
\begin{align*}
\rint S_A&=\{(f_a,f_1,f_b,f_2)\mid(f_1,f_2)\in\rint\Pi(S_A,\F\times\F)\\&\text{ and }(f_a,f_b)\in\rint\Pi(S_A\cap(\F_a\times\pset{f_1}\times\F_b\times\pset{f_2}),\F_a\times\F_b)\}.
\end{align*}
where $S_A=\Pi(\calM,\F_a\times\F\times\F_b\times\F)$. By observing that $\Pi(S_A,\F\times\F)=\Pi(\calM,\F\times\F)=C_f$, we see that the condition \eqref{e:maxmon cond for MA} is equivalent to the existence of $\barf\in\F$ such that $(\barf,-\barf)\in\rint C_f$. Therefore, $\calM_A$ is maximal monotone due to \ref{i:mm cond 1}. As such, it follows from \eqref{e:im L* domain} that $_{B}(\calM_A)$, and thus $\calM_a\rightleftarrows\calM_b$, is maximal monotone if 
\beq\label{e:maxmon cond for BMA}
\im B^*\cap\rint \Pi(\calM_A,\E_a\times\E\times\E_b)\neq\emptyset.
\eeq
To verify this condition, let $S_B=\Pi(\calM_A,\E_a\times\E\times\E_b)$ and note that 
\begin{align*}
\rint S_B&=\set{(e_a,e,e_b)}{e\in\rint\Pi(S_B,\E)\\
&\qquad\text{ and }(e_a,e_b)\in\rint\Pi(S_B\cap(\E_a\times\pset{e}\times\E_b))}.
\end{align*}
Therefore, \eqref{e:maxmon cond for BMA} holds if and only if $0\in\rint\Pi(S_B,\E)$. Note that
\[
\Pi(S_B,\E)=\set{e}{\exists\,(e_1,e_2)\text{ s.t. }e=e_1-e_2}.
\]
As such, \ref{i:mm cond 2} is equivalent to $0\in\rint\Pi(S_B,\E)$ and hence 
\eqref{e:maxmon cond for BMA}. Consequently, $\calM_a\rightleftarrows\calM_b$ is maximal monotone.\EP

Furthermore, maximal {\it cyclical} monotonicity is also preserved under composition as stated in the following theorem.

\bthe\label{t:mcm preserve}
Let $\calM_a \subset \F_a \times \F \times \E_a \times \E$ and $\calM_b \subset \F_b \times \F\times \E_b \times \E$ be maximal cyclically monotone relations that are generated by proper lower semicontinuous convex functions $\phi_a:\F_a\times\F\rightarrow(-\infty,+\infty]$ and $\phi_b:\F_b\times\F\rightarrow(-\infty,+\infty]$, respectively. Let
\[
C_f=\set{(f_1,f_2)}{f_1\in\Pi(\dom\phi_a,\F)\text{ and }f_2\in\Pi(\dom\phi_a,\F)}.
\]
and
\[
C_e=\set{(e_1,e_2)}{\exists\,f\text{ s.t. } (f,e_1)\in\dom\phi_a\text{ and }(-f,e_2)\in\dom\phi_b}.
\]
Suppose that there exists $(\barf,\bare)\in\F\times\E$ such that
\ben[label=(\roman*), ref=(\roman*)]
\item\label{i:mm cond 1-subdiff} $(\barf,-\barf)\in\rint C_f$ and
\item\label{i:mm cond 2-subdiff} $(\bare,\bare)\in\rint C_e$.
\een
Then, $\calM_a \rightleftarrows \calM_b\subset\F_a \times \F_b \times \E_a \times \E_b$ is a maximal cyclically monotone relation that is generated by $\theta^\star:\F_a\times\F_b\rightarrow(-\infty,+\infty]$ where $\theta:\E_a\times\E_b\rightarrow(-\infty,+\infty]$ is given by
\[
\theta(e_a,e_b)=\phi^\star(e_a,0,e_b)
\]
and $\phi:\F_a\times\F\times\F_b\rightarrow(-\infty,+\infty]$ is given by
\[
\phi(f_a,f,f_b)=\phi_a(f_a,f)+\phi_b(f_b,-f).
\]
\ethe
\BP
Let $\calM$, $A$, $B$, $\calM_A$, and $_B(\calM_A)$ be as in the proof of Theorem~\ref{t:mm preserve}. Note that $\calM$ is generated by the proper lower semicontinuous convex function $\phi_{ab}:\F_a \times \F \times \F_b \times \F\rightarrow(-\infty,+\infty]$ given by
\[
\phi_{ab}(f_a,f_1,f_b,f_2)=\phi_{a}(f_a,f_1)+\phi_{b}(f_b,f_2).
\]
For a proper convex function $\Psi:\calG\rightarrow(-\infty,+\infty]$ and a linear map $L:\calH\rightarrow\calG$, let $\psi\circ L:\calH\rightarrow(-\infty,+\infty]$ denote the function given by $h\mapsto \psi(Lh)$.  It follows from the definition of $\calM_A$ that $(f_a,f,f_b,e_a,e,e_b)$ if and only if 
\beq\label{e:char MA-subdiff}
(e_a,e,e_b)\in A^*\partial\phi_{ab}\big(A(f_a,f,f_b)\big).
\eeq
Similar arguments as employed in the proof Theorem~\ref{t:mm preserve} show that \ref{i:mm cond 1-subdiff} is equivalent to
\[
\im A\cap\rint\dom\phi_{ab}\neq \emptyset.
\]
Then, it follows from \cite[Prop. 5.4.5]{bertsekas} that $A^*\partial\phi_{ab}(Ax)=\partial(\phi_{ab}\circ A)(x)$ for all $x\in\F_a\times\F\times\F_b$. Since $\phi_{ab}$ is lower semicontinuous, so is $\phi_{ab}\circ A$. As such, we see from \eqref{e:char MA-subdiff} that $\calM_A$ is maximal cyclically monotone and generated by $\phi=\phi_{ab}\circ A$. Now, it follows from \eqref{e:BMA is Ma comp Mb}, the definition of $_B(\calM_A)$, and \eqref{e:inverse partial duality} that $(f_a,f_b,e_a,e_b)\in\calM_a\rightleftarrows\calM_b$ if and only if
\beq\label{e:char MB-subdiff}
(f_a,f_b)\in B\partial\phi^\star\big(B^*(e_a,e_b)\big)
\eeq
One can show that \ref{i:mm cond 2-subdiff} is equivalent to
\[
\im B^*\cap\rint\dom\phi^\star\neq\emptyset
\]
by employing similar arguments to those in the proof of Theorem~\ref{t:mm preserve}. Then, it follows from \cite[Prop. 5.4.5]{bertsekas} that $B\partial\phi^\star(B^*y)=\partial(\phi^\star\circ B^*)(y)$ for all $y\in\F_a\times\F_b$. Since $\phi^\star$ is lower semicontinuous, so is $\phi^\star\circ B^*$. Consequently, \eqref{e:char MB-subdiff} and \eqref{e:inverse partial duality} imply that $\calM_a \rightleftarrows \calM_b$ is maximal cyclically monotone and generated by $(\phi^\star\circ B^*)^\star$. Since $\theta=\phi^\star\circ B^*$, this concludes the proof.\EP

The following adaptation of Theorem \ref{t:mcm preserve} applies to the {\it alternative} interconnection $e_2=f_3, e_3=f_2$.
Consider two maximal cyclically monotone relations $\M_a, \M_b$ with generating convex functions $g(e_1,e_2)$ and $h(e_3,e_4)$:
\bq
\M_a = \{ (e_1,e_2, f_1=\partial_{e_1}g, f_2=\partial_{e_2}g) \}, \quad \M_b = \{ (e_3,e_4, f_3=\partial_{e_3}h, f_4=\partial_{e_4}h) \}
\eq
Assume $\dim e_2=\dim e_3$, and consider the positive feedback interconnection $e_2=f_3, e_3=f_2$. This yields the relation
\[
\M\! :=\!\{(e_1,e_4,f_1,f_4) \!\mid\! \exists e_2=f_3, e_3=f_2 \mbox{ s.t. }(e_1,e_2,f_1,f_2)\! \in \M_a, (e_3,e_4,f_3,f_4)\! \in \M_b \}.
\]
\begin{proposition}
\label{p:mcm preserve}
The relation $\M$ is maximal cyclically monotone, with generating function
\bq
\label{eq:gener}
(e_1,e_4)\mapsto\inf_{e_2,e_3} \left( g(e_1,e_2) + h(e_3,e_4) - e_2^T e_3 \right)
\eq
\end{proposition}
\BP
The proof can be directly based on the proof of Theorem \ref{t:mcm preserve} formulated for the canonical interconnection $f_2=-f_3, e_2=e_3$. Define $h^*$ as the partial convex conjugate of $h$ with respect to $e_3$, i.e.,
\bq
h^*(f_3,e_4) = \sup_{e_3} \left[e_3^T f_3 - h(e_3,e_4) \right]
\eq
Note that by the definition of the partial convex conjugate
\bq
\partial_{f_3} h^*=e_3, \; \partial_{e_4} h^*= - \partial_{e_4} h
\eq
Hence $\M$ as obtained from $\M_a, \M_b$ via the interconnection equations $e_2=f_3, e_3=f_2$, can be also understood as the canonical interconnection of $\M_a$ and $\M_b^*$, where 
\bq
\M_b^* = \{ (f_3,e_4, e_3=\partial_{f_3}h^*, f_4=  -\partial_{e_4}h^*) \}
\eq
Thus in view of Theorem \ref{t:mcm preserve} $\M$ has the generating function
\bq
\inf_x \left[g(e_1,x) - h^*(x,e_4)\right]
\eq
Substitution of the expression of $h^*$ is immediately seen to result in \eqref{eq:gener}.
\EP

An application of Proposition \ref{p:mcm preserve} is the following. Consider two maximal cyclically monotone port-Hamiltonian systems of the form as given in \eqref{systemconvex}, that is
\begin{equation} 
\begin{array}{rcl}
\dot{x}_i& = & -\frac{\partial K_i}{\partial e_i} (\nabla H_i(x_i),u_i), \quad e_i = \nabla H_i(x_i) \\[2mm]
y_i & = & \frac{\partial K_i}{\partial u_i} (\nabla H_i(x_i),u_i), \qquad i=1,2
\end{array}
\end{equation}
Now interconnect both systems by the positive feedback $u_1=y_2, \, u_2=y_1$. By Proposition \ref{p:mcm preserve} this leads to the maximal cyclically monotone port-Hamiltonian system
\bq
\bma \dot{x}_1 \\[2mm] \dot{x}_2 \ema = - \bma \frac{\partial K}{\partial e_1} \\[2mm]  \frac{\partial K}{\partial e_2} \ema(\nabla H_1(x_1), \nabla H_2(x_2)),
\eq
where the convex function $K(e_1,e_2)$ is given as
\bq
K(e_1,e_2)= \inf_{u_1,u_2} \left( K_1(e_1,u_2) + K_2(e_2,u_2) - u_1^T u_2 \right)
\eq
\begin{example}
Consider two gradient algorithms in continuous time, i.e.,
\beq
\begin{array}{rcl}
\tau_{i} \dot{q}_i  & = & - \frac{\partial P_i}{\partial q_i}(q_i)  - B_iu_i \\[2mm]
y_i & = & B_i^T q_i, \qquad i=1,2
\end{array}
\eeq
which converge for $u_i=0$ to the minimum of the convex functions $P_i(q_i)$. Now consider the {\it coupled} gradient algorithm that is resulting from the interconnection $u_1=y_2, u_2=y_1$. This leads to a maximal cyclically monotone port-Hamiltonian system with respect to the convex function determined as
\bq
\inf_{u_1,u_2} \left(P_1(q_1) + P_1(q_2) + q_1^T B_1u_1 + q_2^T B_2u_2 - u_1^T u_2 \right)
\eq
Clearly the minimum is attained for $u_1=B_1^T q_2 , u_2=B_2^T q_1$, leading to the convex function
\bq
P(q_1,q_2) := P_1(q_1) + P_2(q_2) + q_1^T B_1 B_2^T q_2
\eq
Hence the coupling of the two gradient algorithms computes the minimum of $P(q_1,q_2)$.
\end{example}

\section{When are port-Hamiltonian systems incrementally port-Hamil\-to\-nian, and conversely}
\label{s:ph-inc-ph}
Let us first relate (maximal) (cyclically) monotone relations to the notions of Dirac structures and energy-dissipating relations as used in the definition of port-Hamiltonian systems.

We begin with showing that every Dirac structure is a maximal monotone relation, and maximal cyclically monotone if and only if it belongs to a special subclass of Dirac structures. This special subclass is defined and characterized as follows \cite{hamgraphs}. 
\begin{definition}
\label{d:separable}
A Dirac structure $\mathcal{D}\subset\mathcal{F}\times\mathcal{E}$ is \textit{separable} if 
\begin{equation}\label{separable}
<e_{a}\mid f_{b}>=0\,,\quad \mbox{ for all } (f_{a},e_{a}),(f_{b},e_{b})\in\mathcal{D}
\end{equation}
\end{definition} 
Separable Dirac structures have the following simple geometric characterization \cite{hamgraphs}.
\begin{proposition}
\label{p:separable}
Any separable Dirac structure $\mathcal{D}\subset\mathcal{F}\times \E$ can be written as 
\begin{equation}\label{separable1}
\mathcal{D}=\mathcal{K}\times\mathcal{K}^{\perp}
\end{equation}
for some subspace $\mathcal{K}\subset\mathcal{F}$, where $\mathcal{K}^{\perp}= \{ e \in \E \mid <e \mid f> = 0, \forall f \in \mathcal{K} \}$. Conversely, any
subspace $\mathcal{D}$ as in (\ref{separable1}) for some $\mathcal{K}\subset\mathcal{F}$ is a separable Dirac structure. 
\end{proposition}
\begin{remark}
A typical example of a separable Dirac structure is provided by Kirchhoff's current and voltage laws of an electrical circuit. Indeed, take e.g. $\F$ to be the space of currents, $\mathcal{K}$ the space of currents satisfying Kirchhoff's current laws. Then $\E=\F^*$ is the space of voltages, and $\mathcal{K}^{\perp}$ defines Kirchhoff's voltage laws. Moreover, $<e_{a}\mid f_{b}>=0$ for all $(f_{a},e_{a}),(f_{b},e_{b})\in \mathcal{K}\times\mathcal{K}^{\perp}$ expresses Tellegen's law.
\end{remark}
\begin{proposition}
\label{l:every dirac}
Every Dirac structure $\mathcal{D} \subset \mathcal{F} \times \mathcal{E}$ is maximal monotone. It is maximal cyclically monotone if and only if  $\mathcal{D}$ is separable. If $\mathcal{D}$ is the graph of a mapping $J:\E \to \F$ or $J:\F \to \E$ then $\D$ is cyclically monotone if and only if $J=0$.
\end{proposition}
\BP
Let $\D \subset \mathcal{F} \times \mathcal{E}$ be a Dirac structure. Let $(f_i,e_i)\in \D$ with $i=1,2$. Since $\dualp{e}{f}=0$ for all $(f,e)\in\calD$ due to Remark~\ref{r:dirac}, we obtain by linearity 
\[
\dualp{e_1-e_2}{f_1-f_2}=0.
\]
Therefore, $\calD$ is monotone on $\mathcal{F} \times \mathcal{E}$. Let $\calD'$ be a monotone relation on $\mathcal{F} \times \mathcal{E}$ such that $\calD\subseteq \calD'$. Let $(f',e')\in\calD'$ and $(f,e)\in\calD$. Since $\calD'$ is monotone, $\calD$ is a subspace, and $\calD\subseteq\calD'$, we have
\[0\leq \dualp{e'-\alpha e}{f'-\alpha f}= \dualp{e'}{f'}-\alpha(\dualp{e'}{f}+\dualp{e}{f'})\]
for any $\alpha\in\R$. This means that
$
\dualp{e'}{f}+\dualp{e}{f'}=0,
$
and hence $(f',e')\in\D^{\pperp}=\D$. Therefore, we see that $\D'\subseteq \D$, and thus $\D'=\D$. Consequently, $\D$ is maximal monotone.

Next, let $\calD$ be separable, i.e. $<e_{a}\mid f_{b}>=0$ for all $(f_{a},e_{a}),(f_{b},e_{b})\in\mathcal{D}$. Then it immediately follows from Definition~\ref{d:monotone} that $\calD$ is cyclically monotone. Conversely, let $\calD$ be cyclically monotone. Then take any $(f_i,e_i)\in\calD$ with $i\in\pset{0,1,2}$. It follows from Definition~\ref{d:monotone} that
\[
\dualp{e_0}{f_0-f_1}+\dualp{e_1}{f_1-f_2}+\dualp{e_2}{f_2-f_0}\geq 0.
\]
Since $\dualp{e}{f}=0$ for all $(f,e)\in\calD$ due to Remark~\ref{r:dirac}, we see that
\beq\label{e:3-mon}
\dualp{e_0}{-f_1}+\dualp{e_1}{-f_2}+\dualp{e_2}{-f_0}\geq 0.
\eeq
As $\calD$ is a subspace, $(-f_0,-e_0)\in\calD$. Therefore, we see from \eqref{e:3-mon} that
\[
\dualp{-e_0}{-f_1}+\dualp{e_1}{-f_2}+\dualp{e_2}{f_0}\geq 0.
\]
By summing this inequality and \eqref{e:3-mon}, we obtain $\dualp{e_1}{-f_2}\geq 0$.  By using the fact that $\calD$ is a subspace, we see that $\dualp{e_1}{f_2}=0$, and thus $\calD$ is separable. 

Finally, let $\mathcal{D}$ be the graph of a mapping $J:\E \to \F$. Since $\calD$ is a Dirac structure necessarily $J$ is skew-symmetric. Take again any $(f_i,e_i)\in\calD$ with $i\in\pset{0,1,2}$, where now $f_i=Je_i$. Then if $\calD$ is cyclically monotone
\[
\dualp{e_0}{J(e_0-e_1)}+\dualp{e_1}{J(e_1-e_2)}+\dualp{e_2}{J(e_2-e_0)}\geq 0
\]
Using $\dualp{e_i}{Je_i}=0$ by skew-symmetry of $J$ this yields $\dualp{e_1}{Je_2} \geq 0$ for all $e_1,e_2$, which clearly implies $J=0$. Similarly for $J:\F \to \E$.
\EP

\begin{remark}
As noticed before, a typical example of a separable Dirac structure is given by Kirchhoff's current and voltage laws. In particular, it follows that for any electrical circuit there exists a convex function specifying Kirchhoff's current and voltage laws. Indeed, consider an electrical circuit whose circuit graph is given by the incidence matrix $D$. Identify as above $\F$ with the set of currents $f=I$ through the edges, and $\E=\F^*$ with the set of voltages $e=V$ across the edges. Then Kirchhoff's current laws are given as $DI=0$ and Kirchhoff's voltage laws as $V\in \im D^T$. The convex function generating the resulting separable Dirac structure is given by (see the proof lines of Lemma \ref{l:-r and d})
\beq
\phi (f)=
\begin{cases}
0 & \text{ if }f\in \ker D\\
+\infty & \text{ otherwise}.
\end{cases}
\eeq
\end{remark}

An arbitrary energy-dissipating relation need {\it not} be a (maximal) monotone relation; as was also demonstrated by some of the examples in the previous section. A {\it special} type of energy-dissipating relation that {\it is} a maximal cyclically monotone relation is that of a {\it linear} energy-dissipating relation which is of {\it maximal dimension}. Such an energy-dissipating relation in the port-variables $(f,e)\in\calF\times\calE$ can be represented as a subspace
\begin{equation}\label{reskernel}
\mathcal{R}= \{ (f,e) \in \F \times \E \mid R_f f - R_e e = 0 \},
\end{equation}
where the matrices $R_f, R_e$ satisfy the property
\begin{equation}\label{ressymmetry}
R_f R_e^T = R_eR_f^T \geq 0,
\end{equation}
together with the dimensionality condition
\begin{equation}\label{resdimensionality}
\rank \begin{bmatrix} R_f & R_e\end{bmatrix} = \dim \mathcal{F}.
\end{equation}
First of all, this is seen to define an energy-dissipating relation as follows. By the dimensionality condition (\ref{resdimensionality}) and the equality in (\ref{ressymmetry}) we can equivalently rewrite the kernel representation (\ref{reskernel}) as an image representation
\begin{equation}\label{resimage}
f = R_e^T \lambda, \quad e = R_f^T \lambda.
\end{equation}
That is, any pair $(f,e)$ satisfying (\ref{reskernel}) also satisfies (\ref{resimage}) for some $\lambda$, and conversely, every $(f,e)$ satisfying (\ref{resimage}) for some $\lambda$ also satisfies (\ref{reskernel}). Hence by (\ref{ressymmetry}) for all $(f,e)$ satisfying (\ref{reskernel}) 
\beq
\label{resstructure}
e^Tf =  \left( R_f^T \lambda \right)^T R_e^T \lambda = \lambda^T R_f R_e^T \lambda \geq 0
\eeq
A subspace $\mathcal{R} \subset  \F \times \E$ as in (\ref{reskernel}) where $R_f,R_e$ satisfy \eqref{ressymmetry} and \eqref{resdimensionality} is called a {\it linear resistive structure}. A linear resistive structure can be regarded as a geometric object having properties which are opposite to those of a Dirac structure, in the sense that a Dirac structure can be regarded as the generalization of a {\it skew-symmetric} map, while a linear resistive relation as the generalization of a positive semi-definite {\it symmetric} map. (Geometrically $\mathcal{R}$ defines a Lagrangian subspace of the linear space $\mathcal{F} \times \mathcal{E}$.)

It turns out that every linear resistive structure $\mathcal{R} \subset \mathcal{F} \times \mathcal{E}$ is maximal cyclically monotone. To elaborate further, note that there exists $R=R^T\geq 0$ such that
\beq\label{e:from khatri}
R_eRR_e^T=R_e R_f^T
\eeq
due to \cite[Thm. 2.5]{khatri:76}. In general, $R$ is not unique but the matrix $RR_e^T$ does not depend on the choice of $R=R^T\geq 0$ satisfying \eqref{e:from khatri}. Now, define the extended real-valued convex function
\beq\label{e:conv fun for lrs}
\phi_\calR(f)=
\begin{cases}
\half f^TRf & \text{ if }f\in\im R_e^T\\
+\infty & \text{ otherwise}.
\end{cases}
\eeq
With these preparations, we can state the following characterization for linear resistive structures.

\blem\label{l:-r and d}
Let $\calR\subset\calF\times\calE$ be a linear resistive structure and $\calD\subset\F'\times\calF\times\calE'\times\calE$ be a Dirac structure. Then, the following statements hold:
\ben[label=\normalfont(\alph*), ref=\normalfont(\alph*)]
\item\label{l:-r and d.1} $\calR$ is generated by $\phi_\calR$ and hence is maximal cyclically monotone.
\item\label{l:-r and d.2} The composition $\mathcal{D} \rightleftarrows \mathcal{R}$ is maximal monotone.
\item\label{l:-r and d.3} Any port-Hamiltonian system with Dirac structure and linear resistive structure is maximal monotone port-Hamiltonian.
\een
\elem

\BP
\ref{l:-r and d.1}: Clearly, $\phi_\calR$ is a proper lower semicontinuous function with $\dom\phi=\im R_e^T$ and
\beq\label{e:sub diff for lrs}
\partial\phi_\calR(f)=
\begin{cases}
Rf+\ker R_e & \text{ if }f\in\im R_e^T\\
\emptyset & \text{ otherwise}.
\end{cases}
\eeq
Now, we claim that $\calR$ is generated by $\phi_\calR$ and hence maximal cyclically monotone. To verify this claim, one needs to show that
\beq\label{e:R generated by phi}
\calR=\set{(f,e)}{e\in\partial\phi_\calR(f)}.
\eeq
To see this, first let $(f,e)\in\calR$. Then, we see from \eqref{resimage} that $f=R_e^T\lambda$ and $e=R_f^T\lambda$ for some $\lambda$. Note that 
\beq\label{e:inclusion one}
\im(RR_e^T-R_f^T)\subseteq\ker R_e
\eeq
due to \eqref{e:from khatri}. As such there must exist $\mu\in\ker R_e$ such that $R_f^T\lambda=RR_e^T\lambda+\mu$. Therefore, it follows from \eqref{e:sub diff for lrs} that $e\in\partial\phi_\calR(f)$. This proves that 
\beq\label{e:R generated by phi-1}
\calR\subseteq\set{(f,e)}{e\in\partial\phi_\calR(f)}.
\eeq
To see that the reverse inclusion also holds, let $(f,e)$ be such that $e\in\partial\phi_\calR(f)$. From \eqref{e:sub diff for lrs}, we see that there exist $\lambda$ and $\mu\in\ker R_e$ such that $f=R_e^T\lambda$ and $e=RR_e^T\lambda+\mu$. Since $\ker R_e\subseteq R_f^T\ker R_e^T$ due to \eqref{reskernel} and \eqref{resimage}, it follows from \eqref{e:inclusion one} that $e=R_f^T(\lambda+\theta)$ where $\theta\in\ker R_e^T$. Note that $f=R_e^T\lambda=R_e^T(\lambda+\theta)$. Consequently, we see that 
\beq\label{e:R generated by phi-2}
\set{(f,e)}{e\in\partial\phi_\calR(f)}\subseteq\calR
\eeq
which, together with \eqref{e:R generated by phi-1}, proves \eqref{e:R generated by phi}.

\ref{l:-r and d.2}: Note first that $\calR$ is clearly maximal monotone. Since both $\calD$ and $\calR$ are subspaces, the sets $\Pi(\calD,\calF)$, $\Pi(\calR,\calF)$, $\Pi(\calD,\calF\times\calE)$, and $\Pi(\calR,\calF\times\calE)=\calR$ are all subspaces. As such, the conditions \ref{i:mm cond 1} and \ref{i:mm cond 2} of Theorem~\ref{t:mm preserve} are trivially satisfied by the choices $\barf=0=\bare$. Consequently, the composition $\mathcal{D} \rightleftarrows \mathcal{R}$ is maximal monotone.

\ref{l:-r and d.3}: This immediately follows from the fact that the definition of a port-Hamilto\-nian system entails the composition $\mathcal{D} \rightleftarrows \mathcal{R}$ of $\mathcal{D}$ and $\mathcal{R}$. 
\EP

Furthermore, if we replace in the definition of a port-Hamiltonian system the energy-dissipating relation $\mathcal{R}$ by a relation $\mathcal{R}'$ such that $-\mathcal{R}'$ is monotone, then also $\mathcal{D} \rightleftarrows \mathcal{R}'$ is monotone, and thus we obtain a monotone port-Hamiltonian system (which is however not necessarily port-Hamiltonian).

\section{Steady state analysis of incrementally port-Hamiltonian systems}
\label{s:steady}
In this section we utilize the theory from the previous section to analyze the set of {\it steady states} (for non-zero constant inputs) of an incrementally port-Hamiltonian system and of interconnections of incrementally port-Hamiltonian systems.
For simplicity of exposition, we will denote throughout this section $\mathcal{Y}:= \mathcal{F}_P, \mathcal{U}:= \mathcal{E}_P$, and correspondingly set $y=f_P, u=e_P$.

\subsection{The steady-state input-output relation}
First recall the notion of \\steady-state input-output relation. Consider an input-state-output system $\Sigma$ given as
\begin{equation}\label{systemMEIP}
\Sigma: \quad
\begin{array}{rcl}
\dot{x} & = & f(x,u), \quad x \in \mathbb{R}^n,  u \in \mathbb{R}^m\\[2mm]
y & = & h(x,u), \quad y \in \mathbb{R}^m
\end{array}
\end{equation}
Consider any constant input vector $\bar{u}$ for which there exists an $\bar{x} \in \mathbb{R}^n$ with $0 = f(\bar{x},\bar{u})$, and denote $\bar{y} = h(\bar{x}, \bar{u})$. Then the set of all such pairs $(\bar{y}, \bar{u})$, i.e.,
\beq
\label{ss1}
\mathcal{G}= \{(\bar{y},\bar{u}) \mid \exists \bar{x}, 0 = f(\bar{x},\bar{u}), \bar{y} = h(\bar{x}, \bar{u}) \}
\eeq
is called the {\it steady-state input-output relation} of $\Sigma$. 

In the case of incrementally port-Hamiltonian systems more can be said about the structure of steady-state input-output relations. First we note the following direct applications of Theorem \ref{t:mm preserve} and Theorem \ref{t:mcm preserve} .
\begin{corollary}
Consider an incrementally port-Hamiltonian system with underlying maximal monotone relation $\mathcal{M} \subset \mathcal{F}_S \times \mathcal{E}_S \times \mathcal{Y} \times \mathcal{U}$. Assume $\M$ satisfies
\beq
\label{s1}
0 \in \rint \Pi(\M,\F_S)
\eeq
and
\beq
\label{s2}
\text{ there exists } \bar{e} \text { such that } \bar{e} \in \rint \{ e_S \mid (0,e_S) \in \Pi(\M,\F_S\times\E_S \}
\eeq
Then
\beq
\mathcal{M}_s= \{ (y,u) \mid \exists e_S \mbox{ such that } (0,e_S,y,u) \in \mathcal{M} \}
\eeq
is also a maximal monotone relation. 

Furthermore, in case the maximal monotone relation is cyclically monotone, and thus is given as the subdifferential of some convex function $K(e_S,u)$, then
\beq
\mathcal{M}_s = \mbox{graph } (\partial K_s),
\eeq
where the convex function $K_s: \mathcal{U} \to \mathbb{R}$ is given as
\[
K_s(u)=K^*(0,u)
\]
with $K^*(f_S,u)$ the partial convex conjugate of $K$ with respect to $e_S$.
\end{corollary}

\BP
First note that $\M_s$ is the composition of $\mathcal{M}$ with the trivial maximally monotone relation $\{(f_S=0,e_S) \in \mathcal{F}_S \times \mathcal{E}_S \}$. Thus in order to apply Theorem \ref{t:mm preserve} we need to show that there exists $(\bar{f}_S, \bar{e}_S) \in \F_S \times \E_S$ such that (following the notation of Theorem \ref{t:mm preserve})
\[
\begin{array}{l}
(i) (\bar{f}_S,-\bar{f}_S) \in \rint \D_f  \\[2mm]
(ii) (\bar{e}_S,\bar{e}_S) \in \rint \D_e
\end{array}
\]
where $\D_f = \Pi (\M,\F_S) \times 0$ and $\D_e = \{ (e_1,e_2) \mid (0,e_1) \in \Pi (\M, \F_S \times \E_S \}$. It is easily seen that conditions $(i),(ii)$ reduce to \eqref{s1} and \eqref{s2}.

The rest of the proof follows from Theorem \ref{t:mcm preserve}.
\EP

From now on we will throughout assume that the conditions \eqref{s1}, \eqref{s2} are satisfied, implying that $\M_s$ is maximal monotone.

It is directly seen that the steady-state input-output relation $\G$ of the maximal monotone port-Hamiltonian system with maximal monotone relation $\M$ is {\it contained} in the maximal monotone relation $\M_s$. Indeed, if $\bar{u},\bar{x}, \bar{y}$ is such that $(0,\frac{\partial H}{\partial x}(\bar{x}),\bar{y},\bar{u}) \in \M$ (and thus $(\bar{y},\bar{u}) \in \G$) then clearly $(0,e_S,\bar{y},\bar{u}) \in \M$, where $e_S= \frac{\partial H}{\partial x}(\bar{x})$. Consequently, $\G$ is at least {\it monotone}. 

However $\M_s$ may be {\it larger} than $\G$ since there may not exist for every $e_S$ such that $(0,e_S,\bar{y},\bar{u}) \in \M$ an $\bar{x}$ such that $e_S=\frac{\partial H}{\partial x}(\bar{x})$. In fact, this non-existence of $\bar{x}$ may be due to two reasons. First, the Hamiltonian $H$ may be such that there does not exist for any $e_S$ an $x$ (not necessarily a steady state) such that $e_S=\frac{\partial H}{\partial x}(x)$. Secondly, if such an $x$ exists it may not be a steady state.
A simple example illustrating the second reason is provided by the {\it nonlinear integrator}
\[
\dot{x} = u, \, y = \frac{\partial H}{\partial x}(x),
\]
where we assume that $H$ is a strictly convex function such that the mapping $x \mapsto \frac{\partial H}{\partial x}(x)$ is surjective. This defines an incrementally port-Hamiltonian system with 
\[
\mathcal{M}= \{(f_S,e_S,y,u) \mid f_S=-u,y=e_S \}
\]
Clearly, $\M_s= \{(y,u) \mid u=0 \}$. However the set of steady states $\bar{x}$ is either empty or is given as the {\it singleton} $\{\bar{x} \mid \frac{\partial H}{\partial x}(\bar{x})=0 \}$, implying that also $\G$ is a singleton, and hence not equal to $\M_s$.

\subsection{Equilibrium independent passivity}
Up to now, no conditions were imposed on the Hamiltonian $H$ in the definition of an incrementally port-Hamiltonian system.
If additionally $H$ is {\it strictly convex}, as well as differentiable, then for every $\bar{x}$ the function $S_{\bar{x}}: \mathcal{X} \to \mR$ defined as
\begin{equation}
\label{storagef}
H_{\bar{x}}(x) := H(x) - \frac{\partial H}{\partial x^T}(\bar{x})(x - \bar{x}) - H(\bar{x})
\end{equation}
(as a function of $x$ and $\bar{x}$ also called the Bregman divergence of $H$ \cite{bregman}, or as a function of $x$ alone for fixed $\bar{x}$ the shifted Hamiltonian \cite{passivitybook}) has a strict minimum at $\bar{x}$, and is again strictly convex. Furthermore, 
\[
\frac{\partial H_{\bar{x}}}{\partial x}(x) = \frac{\partial H}{\partial x}(x)- \frac{\partial H}{\partial x}(\bar{x})
\]
Hence for any $(\bar{u},\bar{y})$ in the steady-state input-output relation of an incrementally port-Hamiltonian system one computes
\beq
\label{eq:shifted}
\frac{d}{dt}H_{\bar{x}}=\frac{\partial S_{\bar{x}}}{\partial x^T}(x) \dot{x} = \left( \frac{\partial H}{\partial x^T}(x)- \frac{\partial H}{\partial x^T}(\bar{x})\right)(\dot{x} - 0) \leq (y - \bar{y})^T(u - \bar{u}),
\eeq
implying passivity with respect to the {\it shifted} passivity supply rate $(y - \bar{y})^T(u - \bar{u})$. This was called shifted passivity in \cite{passivitybook}, while the property that this holds for {\it any} steady state values $(\bar{u},\bar{x},\bar{y})$ was coined as {\it equilibrium independent passivity} in \cite{arcak}. Summarizing
\begin{proposition}
Consider a maximal monotone port-Hamiltonian system with respect to the maximal monotone relation
\[
\mathcal{M} \subset  \mR^n \times \mR^n \times \mR^m \times \mR^m,
\]
with a strictly convex differentiable Hamiltonian $H: \mathbb{R}^n \to \mathbb{R}$. Then the system is equilibrium independent passive, with static input-output relation given by the monotone relation $\mathcal{G} \subset \mathcal{G}_{\mathcal{M}}$, and with storage functions $H_{\bar{x}}$ having a strict minimum at $\bar{x}$.
If additionally $H$ is such that for every $\bar{e}_x$ there exists an $\bar{x}$ with $\bar{e} = \frac{\partial H}{\partial x}(\bar{x})$ then $\mathcal{G} = \mathcal{G}_{\mathcal{M}}$.
\end{proposition}
The case $\mathcal{G} = \mathcal{G}_{\mathcal{M}}$ was called {\it maximal equilibrium independent passivity} in \cite{buergerzelazo}. (Maximal) equilibrium independent passivity is a desirable property for showing (asymptotic) stability of the steady state values of a port-Hamiltonian system for different constant input values, since by \eqref{eq:shifted} the shifted Hamiltonians can be employed as Lyapunov functions for $u=\bar{u}$.

\subsection{Determination of the steady state of the interconnection of incrementally port-Hamiltonian systems}
In this subsection we analyze how the steady-state of the interconnection of multiple maximal cyclically monotone port-Hamiltonian systems can be computed, under additional assumptions, by solving a convex {\it optimization} problem. This subsection is motivated by some of the developments in \cite{buergerzelazo}.
  
Consider $k$ maximal monotone port-Hamiltonian systems with input and output vectors $u_i \in \mathbb{R}^{m_i}, y_i \in \mathbb{R}^{m_i}, $ and maximal monotone relations $\M_i, i=1, \cdots,k$. Let, as before, $\mathcal{M}^s_i \subset \mathbb{R}^{m_i} \times  \mathbb{R}^{m_i}, i=1, \cdots,k$ be maximal monotone relations. Additionally, assume that $\mathcal{M}^s_i, i=1, \cdots,k$ are maximal {\it cyclically} monotone, and thus the graphs of subdifferentials of convex functions $K_i(u_i), i=1, \cdots,k$. 

Consider now an interconnection of the following general type. For any subset $\pi \subset \{1, \cdots, k\}$ define
\beq
\begin{array}{rclrl}
f_i &:= & u_i, \quad i \in \pi, \quad f_i &:= & y_i, \quad i \notin \pi \\[2mm]
e_i &:= & y_i, \quad i \in \pi, \quad e_i &:= & u_i, \quad i \notin \pi
\end{array}
\eeq
Furthermore, consider any subspace $\mathcal{C}$ of the linear space of variables $(f_1, \cdots, f_k) \in \mathbb{R}^{m_1} \times \cdots, \mathbb{R}^{m_k}$, and define {\it interconnection constraints}
\beq
(f_1, \cdots, f_k) \in \mathcal{C}, \quad 
(e_1, \cdots, e_k) \in \mathcal{C}^{\perp}
\eeq
The main message of this subsection is that finding a steady state of the interconnected system can be performed by solving a {\it convex minimization problem}. Define the convex function $K: \mathbb{R}^{m_1} \times \cdots \times \mathbb{R}^{m_k} \to \mathbb{R} \cup \{ \infty \}$ given as
\beq
K(f_1,\cdots, f_k):= \sum_{i \in \pi} K_i(u_i) + \sum_{ i \notin \pi}K_i^*(y_i)
\eeq
Now consider the minimization
\beq
\min_{(f_1, \cdots,f_k) \in \mathcal{C}} K(f_1, \cdots, f_k)
\eeq
and write $\mathcal{C} = \ker C$ for some constraint matrix $C= \col (C_1, \cdots,C_k)$. Then the minimization is equivalent to the {\it unconstrained} minimization
\beq
\min_{(f_1, \cdots,f_k), \lambda} K(f_1, \cdots, f_k) - \sum_{i}^s \lambda_i^TC_i f_i
\eeq
where $\lambda$ is a corresponding vector of Lagrange multipliers. This yields the first-order optimality conditions
\begin{equation}
\begin{array}{l}
0 \in \frac{\partial K_i}{\partial u_i}(u_i) - C_i^T \lambda  , \quad i \in \pi \\[2mm]
0 \in \frac{\partial K^*_i}{\partial y_i}(y_i) - C_i^T \lambda  , \quad i \notin \pi
\end{array}
\end{equation}
Consider a solution $(\bar{f}_1, \cdots, \bar{f}_k) \in \mathcal{C}$ of these first-order optimality conditions. Hence there exist $\bar{e}_i= \bar{y}_i \in  \frac{\partial K_i}{\partial u_i}(\bar{u}_i), i \in \pi,$ and $\bar{e}_i= \bar{u}_i \in \frac{\partial K^*_i}{\partial y_i}(\bar{y}_i), i \notin \pi,$ such that $\bar{e} \in \im C^T$, which is nothing else than $\bar{e} \in \mathcal{C}^{\perp}$.

This yields the following theorem regarding the equilibrium of the interconnection of maximal cyclically monotone port-Hamiltonian systems.
\begin{theorem}
Consider $k$ maximal cyclically monotone port-Hamiltonian systems with input and output variables $u_1, \cdots, u_k, y_1, \cdots, y_k$ where $u_i \in \mathbb{R}^{m_i}, y_i \in \mathbb{R}^{m_i}, i=1, \cdots,k$. Assume that the maximal monotone relations ${\mathcal{G}_{\mathcal{M}}}_i \subset \mathbb{R}^{m_i} \times \mathbb{R}^{m_i}$ are given as the graph of subdifferentials $\partial K_i$ for convex functions $K_i,  i=1, \cdots, k$. Furthermore, let $\pi \subset \{1, \cdots, k \}$ be an index set and consider any constraint subspace $\mathcal{C} \subset \mathbb{R}^{m_1} \times \cdots, \mathbb{R}^{m_k}$ leading to the interconnection
\[
(f_1, \cdots, f_k) \in \mathcal{C}, (e_1, \cdots, e_k) \in \mathcal{C}^{\perp}
\]
Then if $(\bar{f}_1, \cdots, \bar{f}_k) \in \mathcal{C}$ is a solution of the minimization 
\[
\min_{(f_1, \cdots,f_k) \in \mathcal{C}} K(f_1, \cdots, f_k)
\]
then there exists $(\bar{e}_1, \cdots, \bar{e}_k) \in \mathcal{C}^{\perp}$. 
\end{theorem}
Note that once we have computed $(\bar{e}_1, \cdots, \bar{e}_k)$ and there exists $(\bar{x}_1, \cdots, \bar{x}_k)$ such that $\frac{\partial H_i}{\partial x_i}(\bar{x}_i)= \bar{e}_i, i=1, \cdots,k$, then this means that $(\bar{x}_1, \cdots, \bar{x}_k)$ is an equilibrium of the interconnected system. Furthermore, if we additionally assume that the Hamiltonians $H_i$ are strictly convex, it follows that this equilibrium is {\it stable}.

Finally, note that the interconnection constraints can be equivalently formulated as the solution of the {\it dual} minimization problem 
\beq
\min_{(e_1, \cdots,e_k) \in \mathcal{C}^{\perp}} K^*(e_1, \cdots, e_k)
\eeq
where 
\beq
K^*(e_1,\cdots, e_k):= \sum_{i \in \pi} K^*_i(y_i) + \sum_{ i \notin \pi}K_i(u_i)
\eeq

\section{Connections with other passivity notions}
\label{s:connections}
In the previous section we already observed that the notion of incrementally port-Hamiltonian systems is closely related to shifted passivity and equilibrium independent passivity. In this section we will discuss how it is closely related to {\it incremental passivity} and {\it differential passivity} as well; at least in case the Hamiltonian is {\it quadratic-affine}. Thus let $H(x)=\half x^TQx + Ax + c$ for some symmetric positive semi-definite matrix $Q$, matrix $A$ and constant $c$. 
In this case, the inequality \eqref{e:disp-like} reduces to
\begin{equation}
\dualp{Q(x_1 - x_2)} {\dot{x}_1-\dot{x}_2} \leq \dualp{e_P^1-e_P^2}{f_P^1-f_P^2}
\end{equation}
which is equivalent to
\begin{equation}
\frac{d}{dt}\half(x_1(t) - x_2(t))^TQ(x_1(t) - x_2(t)) \leq (e_P^1(t)-e_P^2(t))^T(f_P^1(t)-f_P^2(t))
\end{equation}
Recall \cite{Desoer, Pavlov, Angeli} that a system $\dot{x} = f(x,u), y=h(x,u)$ with $x \in \mathbb{R}^n, u,y \in \mathbb{R}^m$ is called {\it incrementally passive} if there exists a nonnegative function $V: \mathbb{R}^n \times \mathbb{R}^n \to \mathbb{R}$ such that
\begin{equation}
\frac{d}{dt} V(x_1,x_2) \leq (u_1 - u_2)^T (y_1 - y_2)
\end{equation}
for all $(x_i,u_i,y_i), i=1,2$ satisfying $\dot{x} = f(x,u), y=h(x,u)$. We immediately obtain the following result.
\begin{proposition}
Any incrementally port-Hamiltonian system with quadratic-affine Hamiltonian $H(x)=\half x^TQx + Ax + c$ with $Q\geq 0$ is {\it incrementally passive}. 
\end{proposition}
\BP
The function $V(x_1,x_2)= \half(x_1 - x_2)^TQ(x_1 - x_2)$ is immediately seen to define an incremental storage function for incremental passivity.
\EP

Recall furthermore from \cite{fsNOLCOS,vdsNOLCOS,fss} the following definition of differential passivity.
\begin{definition}\label{def:differentialpassivity}
Consider a nonlinear control system $\Sigma$ with state space $\X$, affine in the
inputs $u$, and with an equal number of outputs $y$, given as
\begin{equation}\label{eq:system}
  \Sigma : \left.
    \begin{array}{l}
      \displaystyle{\dot{x} = f(x) + \sum_{j=1}^m u_j g_j (x)} \, , \\
      y_j = H_j(x) \, , \quad j=1,\ldots,m \, , 
    \end{array}
  \right.
\end{equation}
The \emph{variational system} along any input-state-output trajectory 
$$t
\in [0,T] \mapsto (x(t),u(t),y(t))$$
is given by the following
time-varying system, cf. \cite{crouch}
\begin{equation}\label{eq:variational}
\begin{array}{rcl}
  \dot{\dx}(t) &=& \pder{f}{x} (x(t)) \dx(t) + \\[2mm]
   && \sum_{j=1}^m u_j(t)
  \pder{g_j}{x}(x(t)) \dx(t) + 
  \sum_{j=1}^m \du_j g_j(x(t)) \\[2mm]
  \dy_j (t) & = & \pder{H_j}{x} (x(t)) \dx(t) \, , \quad j=1,\ldots,m \, ,
  \end{array}
\end{equation}
with state $\dx \in \real^n$, where $\du = (\du_1, \ldots, \du_m)$,
$\dy = (\dy_1,\ldots, \dy_m)$ denote the inputs and the outputs of the
variational system. 
Then $\Sigma$ is called {\it differentially passive} if the system together with all its variational systems is dissipative with respect to the supply rate $\du^T \dy$, that is, if there exists a function $P : T\X \to \mathbb{R}^+$ (called the {\it differential storage function}) satisfying
\begin{equation}\label{eq:difpassive}
\frac{d}{dt} P \leq   \du^T\dy
\end{equation}
for all $x, u, \du$.
\end{definition}
Similar to incremental passivity we obtain
\begin{proposition}
A monotone port-Hamiltonian system with quadratic-affine Hamiltonian $H(x)=\half x^TQx + Ax + c$ with $Q\geq 0$ is differentially passive.
\end{proposition}
\BP
Consider the {\it infinitesimal} version of (\ref{e:disp-like}). In fact, let $(f_P^1, e_P^1,x_1)$ and $(f_P^2, e_P^2,x_2)$ be two triples of system trajectories arbitrarily near each other. Taking the limit we deduce from (\ref{e:disp-like})
\begin{equation}\label{e:disp-like1}
\dx^T \frac{\partial^2 H}{\partial x^2}(x) \delta \dot{x} \leq \delta e_P^T \delta  f_P
\end{equation}
where $\dx$ denotes the variational state, and $\partial f_P, \partial e_P$ the variational inputs and outputs). If the Hamiltonian $H$ is a quadratic function $H(x)=\half x^TQx  + Ax + c$ then the left-hand side of the inequality (\ref{e:disp-like1}) is equal to $\frac{d}{dt} \half \dx^T Q \dx$, and hence amounts to the differential dissipativity inequality
\begin{equation}
\frac{d}{dt} \half \dx^T Q \dx \leq \delta e_P^T \delta  f_P,
\end{equation}
implying that the monotone port-Hamiltonian system is differentially passive, with differential storage function $\half \dx^T Q \dx$.
\EP

Of course, the assumption of a quadratic-affine Hamiltonian $H(x)=\half x^TQx  + Ax + c$ in order to let the monotone port-Hamiltonian system be incrementally passive and differentially passive is restrictive. On the other hand, it is known from the literature \cite{kulkarni,chaffey} that for 'unconditional' incremental properties such an assumption may be necessary as well. For example we can formulate the following simple result. Consider a scalar nonlinear integrator system 
\beq
\label{integrator}
\dot{x}=u, \quad y=\frac{dH}{dx}(x)
\eeq
In order to evaluate its incremental properties consider two copies
\beq
\label{integrators}
\dot{x}_1=u_1, \dot{x}_2=u_2, \quad y_1 = \frac{dH}{dx_1}(x_1), y_2 = \frac{dH}{dx_2}(x_2)
\eeq
Then the system \eqref{integrator} is incrementally passive iff there exists $S(x_1,x_2) \geq 0$ satisfying
\beq
\frac{\partial S}{\partial x_1}(x_1)u_1 + \frac{\partial S}{\partial x_2}(x_2)u_2 \leq \big(u_1-u_2 \big) \big(\frac{dH}{dx_1}(x_1) - \frac{dH}{dx_2}(x_2) \big)
\eeq
for all $x_1,x_2,u_1,u_2$ related by \eqref{integrators}.
This is equivalent to
\beq
\frac{\partial S}{\partial x_1}(x_1,x_2)=\frac{dH}{dx_1}(x_1)-\frac{dH}{dx_2}(x_2) = - \frac{\partial S}{\partial x_2}(x_1,x_2)
\eeq
for all $x_1,x_2$. Differentiation of the first equality with respect to $x_2$, and of the second equality with respect to $x_1$, yields
\beq
- \frac{d^2H}{dx^2_2}(x_2) = \frac{\partial^2 S}{\partial x_1 \partial x_2}(x_1,x_2) = - \frac{d^2H}{dx^2_1}(x_1),
\eeq 
implying that $\frac{d^2H}{dx^2}(x)$ is a constant; i.e., $H(x)$ must be a quadratic-affine function $H(x)= \half qx^2 + ax + c$, for some constants $q,a,c$. Hence the \eqref{integrator} is incrementally passive {\it if and only if} $H$ is quadratic-affine (in which case the integrator is actually linear). This example is easily extendable to more general situations, basically implying that unconditional incremental passivity implies a quadratic-affine storage function.

\section{Conclusions}\label{s:conc}
The notion of an incrementally port-Hamiltonian system was first introduced in \cite{camvds}; basically replacing the composition of a Dirac structure and an energy-dissipation relation in a standard port-Hamiltonian system by a general monotone relation. The present paper discusses the properties of incrementally port-Hamiltonian systems in much more detail; including a wealth of examples and the formulation of specific system subclasses. In particular, the current paper studies the class of maximal {\it cyclically} monotone port-Hamiltonian systems and its connection to convex generating functions. From a mathematical point of view a key contribution of the present paper is a detailed treatment of composition of maximal (cyclically) monotone relations, and its implications for the interconnection of incrementally port-Hamiltonian systems. Indeed, it is shown that under mild technical conditions the composition of maximal (cyclically) monotone relations defines a maximal (cyclically) monotone relation. 

Apart from the abundance of physical examples, this relates incrementally port-Hamiltonian systems to convex optimization as well. Such relations are multi-faceted; from the formulation of gradient and primal-dual gradient algorithms in continuous time as incrementally port-Hamiltonian systems to the computation of the equilibrium of interconnected incrementally port-Hamiltonian systems via convex optimization. Furthermore, apart from the convex generating functions of maximal cyclically monotone relations, {\it another} use of convexity in this incrementally port-Hamiltonian framework is the consideration of convex Hamiltonian functions. The use of the {\it Bregman divergence} of a convex function already turns out to be natural in assessing the stability of steady states of (interconnected) incrementally port-Hamiltonian systems, but much more connections between port-Hamiltonian theory and convex analysis are still to be explored.

The precise dynamical properties of incrementally port-Hamiltonian systems still remain somewhat illusive. The dynamical implications of the key inequality \ref{e:disp-like} are only fully clear if the Hamiltonian $H$ is a quadratic-affine function in suitable coordinates.  Indeed, in this case the incrementally port-Hamiltonian system is incrementally and differentially passive. On the other hand, as shown in \cite{kulkarni} and in the example of a scalar integrator discussed at the end of the previous section, unconditional incremental properties are typically very demanding (see also the theory of contractive systems \cite{lohmiller}), and one could argue that the notion of incrementally port-Hamiltonian systems is less restrictive (although yet less clear from a dynamical perspective).

%

\bibliographystyle{IEEEbib}
\bibliography{inc-PH}

\end{document}